\documentclass[a4paper,11pt]{article}

\usepackage{amsmath}
\usepackage{amssymb}
\usepackage{amsthm}
\usepackage{amsfonts}
\usepackage{mathtools}
\usepackage{graphicx}
\usepackage{tikz,pgfplots}
\usepackage[utf8]{inputenc}
\usepackage{pgfplots}
\usepackage{tikz}
\pgfplotsset{compat=1.9}
\usepackage{subcaption}
\usepackage{caption}

\usepackage{titlesec}
\usepackage{enumitem}
\usepackage{mathtools}
\usepackage{listings}
\usepackage{geometry}
\usepackage{bbm}
\usepackage{xcolor}
\usepackage{hyperref}
\usepackage{caption}
\usepackage{mathrsfs}
\usepackage{booktabs}
\usepackage{floatrow}
\newfloatcommand{capbtabbox}{table}[][\FBwidth]
\usepackage{blindtext}
\usepackage{colortbl}
\definecolor{mycolor4}{RGB}{0,128,128} 
\definecolor{mycolor2}{RGB}{255,127,0} 
\definecolor{darkgreen}{RGB}{0,100,0} 
\definecolor{lightgreen}{rgb}{0.6, 1.0, 0.6}
\definecolor{purple}{rgb}{0.5, 0.0, 0.5}

\newtheorem{theorem}{Theorem}[section]

\newtheorem{proposition}[theorem]{Proposition}

\definecolor{alizarin}{rgb}{0.82, 0.1, 0.26}
\definecolor{atomictangerine}{rgb}{1.0, 0.6, 0.4}
\definecolor{brandeisblue}{rgb}{0.0, 0.44, 1.0}

\newcommand{\re}{\mbox{\rm Re}}

\makeatletter
\newcommand*{\rom}[1]{\expandafter\@slowromancap\romannumeral #1@}
\makeatother

\newcommand{\0}{L^2(\D)}

\newcommand{\3}{H^{1}_{0}(\D)}
\newcommand{\eps}{\varepsilon}

\newcommand{\D}{\mathcal{D}}

\newcommand{\lm}{\lambda}

\newcommand{\R}{\mathbb{R}}
\newcommand{\C}{\mathbb{C}}

\newcommand{\Ltwo}[2]{(#1,#2)_{L^2(\D)}}
\newcommand{\Hone}[2]{(#1,#2)_{H^1(\D)}}

\definecolor{mycolor1}{rgb}{0.00000,0.44700,0.74100}%

\newcommand{\quotes}[1]{``#1''}

\newcommand{\dx}{\hspace{2pt}\mbox{d}x}

\newcommand{\ci}{\mathrm{i}} 

\newcommand{\sR}{\mbox{\rm \tiny A}}

\newcommand{\nablaR}{\nabla_{\hspace{-2pt}\sR}}

\newcommand{\tangentspace}[1]{T_{#1}\mathbb{S}}

\begin{document}

\begin{center}
{\LARGE 
Riemannian conjugate Sobolev gradients and their application to compute ground states of BECs  
\renewcommand{\thefootnote}{\fnsymbol{footnote}}\setcounter{footnote}{0}
 \hspace{-3pt}\footnote{Y. Ai and S. Yuan acknowledge the support by The Chinese University of Hong Kong and P. Henning and M. Yadav acknowledge the support by the German Research Foundation (DFG grant HE 2464/7-1).}}\\[2em]
\end{center}

\begin{center}
{\large Yueshan Ai\footnote[1]{\label{affiliation1}Department of Mathematics, The Chinese University of Hong Kong, Shatin, Hong Kong.\\ email: \href{mailto:ysai@link.cuhk.edu.hk}{ysai@link.cuhk.edu.hk} and \href{mailto:styuan@link.cuhk.edu.hk}{styuan@link.cuhk.edu.hk} }, Patrick Henning\footnote[2]{\label{affiliation2}Department of Mathematics, Ruhr-University Bochum, DE-44801 Bochum, Germany.\\ email: \href{mailto:patrick.henning@rub.de}{patrick.henning@rub.de} and \href{mailto:mahima.yadav@rub.de}{mahima.yadav@rub.de} }, Mahima Yadav\textsuperscript{\ref{affiliation2}} and Sitong Yuan\textsuperscript{\ref{affiliation1}}
}\\[2em]
\end{center}

\begin{center}
{\large{July 6, 2025}}
\end{center}

\begin{abstract}
This work considers the numerical computation of ground states of rotating Bose-Einstein condensates (BECs) which can exhibit a multiscale lattice of quantized vortices. This problem involves the minimization of an energy functional on a Riemannian manifold. For this we apply the framework of nonlinear conjugate gradient methods in combination with the paradigm of Sobolev gradients to investigate different metrics. Here we build on previous work that proposed to enhance the convergence of regular Riemannian gradients methods by an adaptively changing metric that is based on the current energy. In this work, we extend this approach to the branch of Riemannian conjugate gradient (CG) methods and investigate the arising schemes numerically. Special attention is given to the selection of the momentum parameter in search direction and how this affects the performance of the resulting schemes. As known from similar applications, we find that the choice of the momentum parameter plays a critical role, with certain parameters reducing the number of iterations required to achieve a specified tolerance by a significant factor. Besides the influence of the momentum parameters, we also investigate how the methods with adaptive metric compare to the corresponding realizations with a standard $H^1_0$-metric. As one of our main findings, the results of the numerical experiments show that the Riemannian CG method with the proposed adaptive metric along with a Polak--Ribi\'ere or Hestenes--Stiefel-type momentum parameter show the best performance and highest robustness compared to the other CG methods that were part of our numerical study.
\end{abstract}


\section{Introduction}
An extreme state of matter with remarkable superfluid properties is formed when a dilute bosonic gas condenses at temperatures close to 0 Kelvin to a so-called Bose--Einstein condensate (BEC), cf. \cite{Bos24,Ein24,PiS03}. Its extraordinary superfluid nature can be checked by verifying the existence of quantized vortices in the rotating BEC. In practical setups, the appearance of such vortices is crucially related to the interplay of a (magnetic or optical) trapping potential $V$ (to confine the condensate) and the angular frequency $\Omega$ of a stirring potential (to rotate the condensate). If the angular frequency is too small compared to the strength of the trapping potential, no vortices appear. If the angular frequency is too high, the condensate can be destroyed by centrifugal forces. Only in an intermediate regime, a rich landscape of vortex pattern can observed and studied. In this paper, we are concerned with the numerical computation of such pattern by seeking ground states (i.e. the lowest energy states) of BECs in a rotating frame.

On a given computational domain $\D \subset \R^d$ (for $d=2,3$) a ground state is mathematically described through its quantum state $u : \D \rightarrow \mathbb{C}$, whereas the vortices become visible in the corresponding density $|u|^2 : \D \rightarrow \mathbb{R}$. The density is usually normalized such that the total mass of the BEC fulfills $\int_{\D} |u|^2 \, \mbox{d}x  =1$ or, more precisely, it should hold
\begin{align*}
u \, \in \,\mathbb{S} := \{  \, v \in H^1_0(\D,\mathbb{C}) \, |  \,\, \| v \|_{L^2(\D)} = 1  \, \},
\end{align*}
where $H^1_0(\D,\mathbb{C})$ denotes the usual Sobolev space of weakly differentiable, square-integratable and complex-valued functions with a vanishing trace $v_{\vert \partial \D}=0$.
In a given configuration, a corresponding ground state is characterized as a global minimizer of the total energy of the system. This energy is described by the Gross--Pitaevskii energy functional $E : \mathbb{S} \rightarrow \mathbb{R}$ given by
\begin{align}
\label{energy1}
E(u)
:= \frac{1}{2}\int_{\mathcal{D}} |\nabla u|^2 + V\, |u|^2  - \Omega\, \bar{u}\, \mathcal{L}_{3} u + \frac{\kappa}{2} |u|^4 \dx
\end{align} 
for $u\in \mathbb{S}$.
Here, $\overline{u}$ denotes the complex conjugate of $u$, $V$ represents the external trapping potential, $\mathcal{L}_3 := - \ci \left( x_1 \partial_{x_2} - x_2 \partial_{x_1} \right)$ denotes the $x_3$-component of the angular momentum (hence describing a rotation around the $x_3$-axis), $\Omega \in \mathbb{R}$ is the corresponding angular velocity of the stirring potential and $\kappa \in \mathbb{R}^+$ is a repulsion parameter that encodes the strength of particle interactions. For a comprehensive introduction to the basic theory and mathematical properties of ground states of rotating BECs, we refer to the papers by Bao et al. \cite{Bao14, BaC13b, BWM05}.

In a nutshell, the finding of a ground state requires to find a minimizer of the energy $E$ on the Riemannian manifold $\mathbb{S}$ (which we recall as incorporating the mass normalization constraint for $u$). Hence, we are concerned with a Riemannian minimization/optimization problem, which is precisely the perspective that we are taking in this paper. Alternatively, the problem can be also viewed through the lens of nonlinear eigenvalue problems by examining the corresponding Euler-Lagrange equations, known as Gross-Pitaevskii eigenvalue problem, associated with the constrained energy minimization problem \cite{AHP21NumMath,JarKM14,PH24,Can00,CaL00,CaL00B,DiC07}. The links between the two perspectives (nonlinear eigenvalue problem vs. Riemannian minimization problem) are elaborated in the review paper \cite{HenJar24}. We also note that the Euler--Lagrange equations can be tackled with Newton-type methods \cite{APS23Newton,WWB17,XXXY21}.

As mentioned above, we adopt in the following the perspective of directly minimizing the energy in an iterative process through discretized Riemannian gradient flows / Riemannian gradient methods, cf. \cite{AltPetSty22,ALT17,BaD04,BWM05,CDLX23,CLLZ24,DaK10,DaP17,KaE10,Zhang2022}. The development of optimization techniques on Riemannian manifolds goes back to Smith et al. \cite{smith2014,Smith_MIT_1998} and has been significantly extended in the past three decades. Modern methods combine the concepts of Sobolev gradients, Riemannian descent directions, Riemannian vector transport, and retractions to the constrained manifold $\mathbb{S}$, cf. \cite{DaP17,TCWW20}. In this paper, we will further explore this path by constructing new {\it metric-adaptive Riemannian conjugate gradient} (CG) {\it methods} for the considered application of rotating BECs. In particular, we will numerically investigate the performance of the new methods. Our experimental analysis primarily focuses on improving the computational efficiency of the schemes, specifically by (a) selecting appropriate metrics for the constrained manifold (or more precisely, the tangent spaces at the current iterates $u^n$) and (b) choosing the proper momentum parameter (denoted as $\beta$) to update the next search direction. With this, we also want to examine an open question posed in \cite{DaP17}: {\it ``it remains an open question whether the use of a well-adapted Riemannian metric defined on the constraint manifold could further improve the performance of the approach.''} 
Note here that the choice of the metric is a crucial ingredient, because it changes the Riemannian gradient of $E$ on the manifold and hence the way how a steepest descent/ascent is characterised. More precisely, the Riemannian gradient of $E$ in a point $v\in \mathbb{S}$ with respect to an $X$-metric (induced by an inner product $(\cdot,\cdot)_X$) is obtained in two steps: First, construct the Riesz representation of $E^{\prime}(v)$ in $X$. The representation is called a Sobolev gradient of $E$ (cf. \cite{Neu97}) and is denoted by $\nabla_X E(v)$. Second, project $\nabla_X E(v)$ into the tangent space at $v \in \mathbb{S}$ with the $X$-orthogonal projection. The resulting function is the Riemannian gradient of $E$ in $v$ w.r.t. to the $X$-metric. In the context of rotating BECs, the most popular metrics for conjugate Riemannian gradient methods are the $L^2$-metric (in combination with suitable preconditioners) as studied by Tang et al. \cite{ALT17,ShuTangZhangZhang2024} and an energy-metric (based on an inner product of the form 
$(u,v)_{\0} + (\nablaR u , \nablaR v )_{\0}$ with $\nablaR=\nabla v - \ci\, \Omega \, \mathbf{A}^{\hspace{-2pt}\top}$ and $\mathbf{A}=(x_2,-x_1,0)$) proposed by Danaila et al. \cite{DaK10,DaP17}. The idea of using adaptively changing metrics for (general) optimization problems is discussed by Ring and Wirth \cite{RingWirth12} and Sato \cite{Sato_2022}. In this work, we investigate a particular adaptively changing metric (based on the location on $\mathbb{S}$) that is selected in such a way that the corresponding Sobolev gradient fulfills $\nabla_{X}E(v) =v$ for all $v \in \3$ in the sense of optimal preconditioning. Motivated by guaranteed energy dissipation and global convergence for the corresponding gradient flows, this adaptive choice for the metric was first suggested in \cite{HeP20} for the Gross--Pitaevskii energy without rotation and later transferred to the case with rotation in \cite{PHMY24_1} and spinor BECs in \cite{WuLiuCai2025}. However, the previous works are only concerned with regular Riemannian gradient methods and their analysis and extensions to conjugate gradient versions were left open. Even though the generalization of these concepts to Riemannian conjugate gradients was recently discussed in \cite{HenJar24}, it was only for $\Omega=0$ and for a Fletcher--Reeves-type momentum parameter that we found to perform suboptimal in our experiments with rotating BECs. 

In this work, we hence close this gap and formulate a corresponding Riemannian conjugate gradient method with adaptively changing metric as sketched above and with various choices for momentum parameters. Selecting a suitable momentum parameter is of great practical importance since different choices can result in significant variations in the number of iterations required to attain a specified tolerance. We will consider four different choices according to the most popular types of momentum parameters, which are Dai--Yuan \cite{dai1999}, Fletcher--Reeves \cite{fletcher1964}, Hestenes--Stiefel \cite{hestenes1952} and Polak--Ribi\'{e}re \cite{polak1969}. If the momentum parameter is uniformly zero, the method reduces to a Riemannian gradient method. Besides comparing the performance of the parameters, we also compare the iteration numbers for the adaptive metric with the iteration numbers for the $H^1_0$-metric. In our numerical experiments we find that the adaptive metric leads to a significant acceleration compared to the $H^1_0$-method and that the corresponding Polak--Ribi\'ere and Hestenes--Stiefel-type parameters performed  significantly better than the other two (including the generic choice $\beta=0$).

So far, the mathematical convergence analysis of Riemannian conjugate gradient methods for the Gross--Pitaevskii problem is fully open and hence beyond the scope of this paper. However, in the case $\beta=0$ (i.e. the case of Riemannian gradient methods), a convergence analysis was recently established in \cite{PHMY24_1}.

Finally we note that a practical application of Riemannian gradient methods also requires a space discretization (of the partial differential equations that have to be solved in each iteration of the gradient method). Typical approaches include spectral and pseudo-spectral methods \cite{AnD14,Bao14,BaC13b,CCM10}, finite element methods \cite{CGHZ11,DaH10,EngGianGrub22,HLP24,HSW21,HeP23,PHMY24} or spectral element methods \cite{CLLZ24-discrete,CLLZ24-discrete}. We leave the choice open in this paper and only investigate the outer gradient iteration. The algorithms can be afterwards combined with any preferred spatial discretization that potentially exploits the structure of the considered metric.

In fact, exploiting the metric to construct spatial approximation spaces can be a powerful tool that is well-established in the treatment of general multiscale problems. In this context, a conventional approximation space, e.g. a finite element space $V_h$, is represented in a different metric in order to improve its approximation properties. This amounts to applying a differential operator to $V_h$ and taking the image space as a new approximation space. For example, if $\kappa=0$ in \eqref{energy1}, then $\mathcal{A}v:=E^{\prime}(v)$ is a linear operator and the Sobolev gradient in the $\mathcal{A}$-metric becomes $\nabla_X E(v):= \mathcal{A}^{-1} E^{\prime}(v)$. In this case, the representation of $V_h$ in the $\mathcal{A}$-metric is given by $\mathcal{A}^{-1}V_h$. For linear problems, this approximation space is well-known as an LOD (Localized Orthogonal Decomposition) space, cf. \cite{HauckPeterseim23}, and its approximation properties for multiscale problems were extensively studied in the literature, see \cite{ActaLOD21,MaP14,LODbook21} and the references therein. In particular, LOD multiscale spaces have been successfully applied to the simulation of Bose--Einstein condensates \cite{HeMaPe14,HePer20,Johan22,Johan24,Johan24b}. Besides the choice $\mathcal{A}^{-1}V_h$, there are other ways to enrich $V_h$ by multiscale features through a suitable metric/operator. Exemplarily, we mention the multiscale finite element method (MsFEM) \cite{HW97,EH09}, the generalized MsFEM (GMsFEM) \cite{ChEfHo23,EfGaHo13}, the multiscale spectral generalized finite element method (MS-GFEM) \cite{BaLi11,BaLiSiSt20,MaSc22} or the wavelet-based edge multiscale finite element methods (WEMsFEM) \cite{FuChLi19,FuChLi25,FuLiCrGu21}.
In the light of this perspective, the usage of metric-driven gradient methods (as discussed in this paper) can be seen as a complementary new path to the well-known concept of metric-driven approximation spaces that are commonly used to solve multiscale problems.

$\\$
\textbf{Outline:} The paper is organised as follows. The model framework is introduced in Section \ref{section:assumption_model_problem}. In Section \ref{section:RCGM}, we recall the concepts of Sobolev gradients and Riemannian steepest descent.  With this we propose a class of Riemannian conjugate gradient methods with different realizations based on the choice of the metric and the momentum parameter. Lastly, in Section \ref{section:num_exp}, we evaluate the performance of the aforementioned Riemannian CG methods for different model problems and discuss the numerical observations.

\section{Setting and problem formulation}
\label{section:assumption_model_problem}
In this section, we introduce the basic notation and state the Gross--Pitaevskii equation for rotating Bose--Einstein condensates, including the necessary assumptions for the existence of ground states. These assumptions are valid for the whole manuscript.
\subsection{Notation and assumptions} 
For the subsequent discussion on minimizing the energy $E$ given by \eqref{energy1}, we consider a bounded Lipschitz-domain $\D \subset \mathbb{R}^d$ for $d=2,3$. Furthermore, we make the following set of assumptions.
\begin{enumerate}[label={(A\arabic*)}]
\item\label{A1} The repulsion parameter  $\kappa \ge 0$ that characterizes the particle interactions is a real-valued constant. The trapping potential $V$ is non-negative, real-valued and essentially bounded on $\D$, i.e.,  $V \in L^{\infty}(\mathcal{D},\mathbb{R}_{\ge 0} )$. To ensure that the centrifugal forces to do not exceed the strength of the trapping potential,  the angular velocity $\Omega \in \R$ must be small enough such that there is a constant $\eps > 0$ with
\begin{align*}
	V(x) - \frac{1 + \varepsilon}{4} \Omega^2 (x_1^2 + x_2^2) \ge  0 \quad \text{for almost all } x \in \D.
\end{align*}
Note here that only the trapping frequencies in the $(x_1,x_2)$-plane are relevant, because the rotation of the BEC is around the $x_3$-axis.
\end{enumerate}
The above assumptions are not only needed to ensure the existence of minimizers (ground states) but also to ensure well-posedness of the adaptive metric that we construct later on. Regarding the existence of ground states, we note that the positivity of $V$ is not required, but only introduced to avoid technical issues in the presentation of our results. Note that assuming positivity of $V$ is not restrictive since constant shifts of the energy do not change the minimizers. Hence, we can just add a (sufficiently large) constant to $V$ to make it positive. The assumption of positivity of $\kappa$ can be relaxed to some extend to small negative values. However, this regime of attractive particle interactions is still not yet fully understood, which is why we exclude it here. Finally, the balancing assumption of $V$ and $\Omega$ is crucial and there exist typically no ground states if it is violated, cf. \cite{BWM05}.\\[-0.5em]

As we are concerned with minimizing a real-valued functional $E$, it is important to note that we must consider $\R$-Hilbert spaces to ensure Fr\'echet-differentiability of $E$ (cf. \cite{AHP21NumMath,Begout22} for a detailed motivation). For that reason, we equip the Sobolev space $H_{0}^{1}(\mathcal{D}):= H_{0}^{1}(\mathcal{D},\mathbb{C})$, over which we minimize $E$, with the real inner product $\Hone{v}{w}:= \re (  \int_{\mathcal{D}} \nabla v \cdot \overline{\nabla w} \dx )$. Recall here $\overline{w}$ as the complex conjugate of $w$.  Similarly, we define the Lebesgue space $L^2(\mathcal{D}):= L^2(\mathcal{D},\C)$ as a real Hilbert space with inner product $\Ltwo{v}{w}:= \re ( \int_{\D} v \, \overline{w} \dx )$. The corresponding $(\text{real})$ dual space is denoted by $H^{-1}(\mathcal{D}) := \big(H_{0}^{1}(\D)\big)^{*}$ with canonical duality pairing $ \langle \cdot , \cdot \rangle := \langle \cdot ,\cdot \rangle _{H^{-1}(\mathcal{D}) , H^{1}_{0}(\mathcal{D})}$.  
%
%

\subsection{Model Problem and first and second order optimality conditions}
\label{Model_setting}
With the notation above, we recall the Gross--Pitaevskii energy functional $E :H_{0}^{1}(\D) \rightarrow \mathbb{R}$ from \eqref{energy1} as
\begin{align*}
E(v)
= \frac{1}{2}\int_{\mathcal{D}} |\nabla v|^2 + V\, |v|^2  - \Omega\, \bar{v}\, \mathcal{L}_{3}v + \frac{\kappa}{2} |v|^4 \dx
\end{align*}
and we recall that a {\it ground state} $u$ is defined as a {\it global minimizer} of $E$ on the restricted $L^2$-unit sphere $\mathbb{S} = \{ v \in H^1_0(\D) \,| \, \|v\|_{L^2(\D)} = 1 \}$. Since $\mathbb{S}$ is a Riemannian manifold, we are hence concerned with a Riemannian optimization problem which reads:
\begin{align}
\label{minimization-problem}
\text{Find  } \hspace{2mm} u \in \mathbb{S} \hspace{2mm} \text{such that } \hspace{2mm} E(u) = \underset{v \in \mathbb{S}} {\text{ min }} E(v).
\end{align}
The homogeneous Dirichlet boundary conditions imposed for $u$ through $H^1_0(\D)$ can be practically justified with the exponential decay of ground states outside of any sufficiently large domain (whose necessary diameter can be estimated through the Thomas--Fermi radius of the condensate), cf. \cite{BaC13b}.

As stated in the introduction, ground states are the stable stationary states at the lowest possible energy level, with a constraint on the mass of the condensate, represented by the condition $u  \in \mathbb{S}$. Local minimizers or saddle points of $E$ on $\mathbb{S}$ that are no ground states are called excited states. Under the assumptions stated in \ref{A1}, the energy functional is positively bounded from below on $\mathbb{S}$ and the existence of global minimizers can be established with standard arguments. In particular, we have the following existence result (cf. \cite{BWM05} and \cite{PHMY24}).
\begin{proposition}[Existence of ground states]
Let $\mathcal{D} \subset \mathbb{R}^d, \, d=2,3$ be a bounded Lipschitz-domain, and assume \ref{A1} holds. Then there exists at least one ground state $u \in \mathbb{S}$ to problem \eqref{minimization-problem} and it holds $E(u)>0$.
\end{proposition}
Since the energy $E$ is invariant under complex phase shifts, i.e., $E(u) = E(e^{\ci \omega} u)$ for all $\omega \in (0, 2\pi]$, we have that $e^{\ci \omega} u$ is a ground state if $u$ is a ground state. Hence, we cannot hope for uniqueness of minimizers in \eqref{minimization-problem}. Although the density $|u|^2$ is independent of such phase shifts, even uniqueness of $|u|^2$ can only be guaranteed up to some (small) critical frequency of $\Omega$, cf. \cite[Section 6.2]{BWM05}. 

In order to set the stage for the Riemannian gradient descent in the next section, we need to consider the derivatives of $E$. It is easy to see that the energy is infinitely often ($\R$-)Fr\'echet differentiable on the Hilbert space $\3$ and a calculation of the first two derivatives yields
\begin{align}
\label{defEprime} \langle E^{\prime}(u) , v \rangle
&= ( \nabla u , \nabla v)_{L^2(\D)} + ( V\, u  - \Omega\, \mathcal{L}_{3} u  + \kappa  |u|^2 u , v )_{L^2(\D)} \, ,  \\
 \langle E''(u) v, w \rangle 
 \label{defEprimeprime}&=  ( \nabla v , \nabla w)_{L^2(\D)} + ( V\, v - \Omega \, \mathcal{L}_{3} v +\kappa \, |u|^2 v  , w )_{L^2(\D)}  +  2\, \kappa \, ( \re (u \overline{v}) \, u,w )_{\0} .
\end{align}
for any $u,v,w \in \3$. Note that $E''(u)$ is a self-adjoint operator with real and positive eigenvalues if \ref{A1} is fulfilled, cf. \cite{PHMY24}.\\[0.3em]
\noindent 
{\bf First order optimality condition.}
Since a ground state $u \in \mathbb{S}$ is a constrained minimizer of $E$ (with constraint $\int_{\D} |u|^2 \dx =1$), we can also formulate the corresponding Euler Lagrange equations. To precise, for any ground state $u \in H^1_0(\D)$, there exists a Lagrange multiplier $\lambda >0$ such that
\begin{align}
\label{eigenvalue-E-derivative-form}
\langle E^{\prime}(u) , v \rangle = \lambda \, ( u , v )_{L^2(\D)}  \qquad \mbox{for all } v\in H^1_0(\D).
\end{align}
The right hand side of \eqref{eigenvalue-E-derivative-form} should be seen as the $\R$-Fre\'chet derivative of  the constraint functional $\tfrac{1}{2}( \| v \|_{L^2(\D)}^2 - 1)$. Since $\lambda$ can be equivalently interpreted as an eigenvalue of $E^{\prime}$, problem \eqref{eigenvalue-E-derivative-form} is commonly known as the {\it Gross--Pitaevskii eigenvalue problem} (GPEVP). By using the explicit representation of $E^{\prime}$ and expressing it in the sense of distributions, we can write the GPEVP \eqref{eigenvalue-E-derivative-form} as: Find $u \in H^1_0(\D)$ and $\lambda \in \R$ such that 
\begin{align}
\label{eigen_value_problem_1}
-\Delta u + V\, u - \Omega \, \mathcal{L}_{3}u + \kappa \, |u|^2u = \lambda \, u.
\end{align}
The eigenvalue $\lambda$ corresponding to a ground state $u$ is called a {\it ground state eigenvalue}. Unlike for non-rotating condensates, i.e., $\Omega=0$, the ground state eigenvalue $\lambda$ may not be the smallest eigenvalue if the angular velocity is sufficiently high. Numerical evidence supporting this were first pointed out in  \cite[Section 6.3]{AHP21NumMath} and can be also found in Section \ref{section:num_exp} of this paper.

Note that \eqref{eigenvalue-E-derivative-form} (and equivalently \eqref{eigen_value_problem_1}) is nothing but the {\it first order optimality condition} for minimizers of $E$ on $\mathbb{S}$. In the case of $\Omega=0$, this condition is even sufficient to characterize global minimizers if $\lambda$ is the smallest eigenvalue of $E^{\prime}$, cf. \cite{CCM10}. For the case $\Omega \not =0$, this is unfortunately no longer possible since $\lambda$ can appear anywhere in the spectrum of $E^{\prime}$. However, we can consider the {\it second order optimality conditions} for our minimization problem.\\[0.3em]
\noindent 
{\bf Second order optimality condition.} For the second order optimality condition, we consider the tangent space at $u$ on the manifold $\mathbb{S}$ which is given by 
$$
\tangentspace{u} := \{ \, v \in H^1_0(\D) \, | \, (u,v)_{\0} = 0 \, \}.
$$
The above tangent space is obtained as the null space of the Fr\'echet derivative of the mass constraint functional $u \mapsto \big( \int_{\D}|u|^2 \dx -1 \big)$ on $\3$. To identify if a given $u \in \mathbb{S}$ is a minimizer, we have to study the spectrum of $E^{\prime\prime}(u)\vert_{\tangentspace{u}}$. To be precise, we seek eigenfunctions $v_i \in \tangentspace{u}$ and corresponding eigenvalues $\lambda_i \in \R$ such that
\begin{align}
\label{eigenvalue-problem-secE}
\langle E''(u) \, v_i , w \rangle = \lambda_i \, (v_i ,w)_{L^2(\D)}
\qquad \mbox{for all } w\in \tangentspace{u}.
\end{align}
Let the eigenvalues be ordered by size with $0<\lambda_1 \le \lambda_2 \le \dots$. If $u$ is a local minimizer of $E$ on $\mathbb{S}$ with Lagrange multiplier $\lambda$, then the {\it necessary second order optimality condition} states that it must hold 
\begin{align*}
 \lambda_i \ge \lambda \quad  \mbox{ for all } i \ge 1, 
\end{align*}
whereas the {\it sufficient second order optimality condition} demands
\begin{align*}
\lambda_1 = \lambda \qquad \mbox{and} \qquad
 \lambda_i > \lambda \quad  \mbox{ for all } i \ge 2.
\end{align*}
Here, $\lambda$ denotes the ground state eigenvalue to the ground state $u$ in the sense of the first order condition (GPEVP) \eqref{eigenvalue-E-derivative-form}. Note that $\lambda$ must always appear at the bottom of the spectrum (if $u$ is a minimizer) since we have $\lambda_1=\lambda$ for the eigenfunction $v_1 = \ci \, u \in \tangentspace{u}$, cf. \cite{PHMY24}. This eigenfunction originates from the aforementioned invariance of $E$ under phase shifts of $u$, i.e. $E(u) = E(e^{\ci \omega} u)$. Consequently, if we have an eigenfunction $u \in \mathbb{S}$ of $E^{\prime}$ with eigenvalue $\lambda$ (i.e. it holds \eqref{eigenvalue-E-derivative-form}) and if we find that $\lambda_1=\lambda$ is a simple eigenvalue of $E^{\prime\prime}(u)\vert_{\tangentspace{u}}$, then the {\it sufficient second order optimality condition} holds and $u$ must be a local minimizer. We call such minimizers quasi-isolated and it is an open conjecture if all minimizers are quasi-isolated.\\[-0.4em]

To give a brief motivation for the first and second order optimality conditions, we can consider arbitrary smooth curves $\gamma(t)$ on $\mathbb{S}$. If, for example $\gamma : (-1,1) \rightarrow \mathbb{S}$ is a two-times differentiable curve with $\gamma(0)=u$, then it must hold $\gamma^{\prime}(0) \in \tangentspace{u}$, since the curve is tangential to $u$ in $t=0$. Furthermore, the map $t \mapsto E(\,\gamma(t)\,)$ must have a minimum in $t=0$. Consequently, it holds
 \begin{align*}
\frac{ \mbox{d} }{ \mbox{d}\hspace{1pt}t} E(\gamma(t)) |_{t=0} = 0 \qquad \mbox{and} \qquad \frac{ \mbox{d}^{\,2} }{ \mbox{d}\hspace{1pt}t^{2}} E(\gamma(t)) |_{t=0} \ge 0.
\end{align*}
By computing the corresponding derivatives and defining $v:=\gamma^{\prime}(0) \in \tangentspace{u}$ (which can be arbitrary), we obtain the first order optimality condition $\frac{ \mbox{d} }{ \mbox{d}\hspace{1pt}t} E(\gamma(t)) |_{t=0} = 0$ as $\langle E^{\prime}(u) , v\rangle =0$ for all $v\in \tangentspace{u}$. Equivalently, we can decompose an arbitrary $w \in H^1_0(\D)$ into $w = \alpha u + v $ for some $\alpha \in \R$ and $v\in \tangentspace{u}$. By defining $\lambda:=\langle E^{\prime}(u) , u \rangle$ we obtain from $\langle E^{\prime}(u) , v\rangle =0$ that
\begin{align*}
\langle E^{\prime}(u) , w\rangle =  \langle E^{\prime}(u) , \alpha u \rangle = \alpha \, \lambda \overset{u \in \mathbb{S}}{=}  \lambda \, ( u ,  \alpha u )_{L^2(\D)} \overset{v\in \tangentspace{u}}{=} \lambda \, ( u , w )_{L^2(\D)},
\end{align*}
which is precisely the GPEVP \eqref{eigenvalue-E-derivative-form}. In a similar fashion, the second order optimality condition is obtained as $\langle E^{\prime\prime}(u) v , v\rangle \ge \, \lambda \, (v,v)_{L^2(\D)}$ for all $v\in \tangentspace{u}$ from the condition $\frac{ \mbox{d}^{\,2} }{ \mbox{d}\hspace{1pt}t^{2}} E(\gamma(t)) |_{t=0} \ge 0$. We refer to \cite{PHMY24,PHMY24_1} for further details for these standard calculations.

In Section \ref{section:num_exp}, we verify that our numerically computed approximations are indeed local minimizers of $E$ by computing the corresponding ground state eigenvalue $\lambda$ and verifying afterwards the sufficient second order optimality condition: $\lm_1=\lm$ and $\lambda_i > \lambda$ for all $i \geq 2$, where $\lm_i$ are the eigenvalues of $E^{\prime\prime}(u)\vert_{\tangentspace{u}}$.

\section{Riemannian conjugate Sobolev gradient methods}
\label{section:RCGM}

In this section we introduce the Riemannian conjugate Sobolev gradient method to minimize the Gross-Pitaeveskii energy functional in a rotating frame. 

\subsection{Riemannian $X$-Sobolev gradients}
The usage of Sobolev gradients gives us a tool to select different metrics for the gradient method in a natural way. In the following, we shall denote by $X$ a Hilbert space equipped with an inner product $(\cdot , \cdot)_{X}$ and with $\3 \subset X$ in the sense of continuous embeddings. The space $X$ gives rise to the corresponding Sobolev gradient $\nabla_X E(v)$ of $E$ in $v$ and we will investigate the influence of different choices for $X$ on the computational efficiency of the method. We also explore different forms of this gradient method based on the choice of the momentum parameter (denoted as $\beta$), which influences the next search direction in our setting.

The basic idea of a Riemannian gradient method is to move, from a given current point on the manifold, into the direction of the Riemannian steepest descent. The steepest descent is given by the negative Riemannian gradient of $E$ in that given point and depends on the metric through the chosen Sobolev gradient. To make this dependency clear, we will call it the {\it Riemannian Sobolev gradient} of $E$. In a point $u \in \mathbb{S}$ on the manifold, the Riemannian Sobolev gradient is given by $P_{u,X} (\,\nabla_X E(u) \,)$ and is obtained in two steps as follows. First, we construct the ($X$-)Sobolev gradient $\nabla_X E(u)$ of $v \in \3$ as the Riesz-representation of the \quotes{regular} gradient $E'(v) \in H^{-1}(\D)$ in $X$. To be precise, we let $\nabla_X E(u) \in X$ denote the unique function with 
\begin{align}
\label{def-Sobolev-gradient}
(\nabla_X E(u) , v )_X = \langle E^{\prime}(u) , v \rangle \qquad \mbox{ for all } v \in \3.
\end{align}
Note that this problem is in general only well-posed if $X$ contains the elements of $H^1_0(\D)$, so that $E^{\prime}(u) \in X^{\ast}$. However, this can be often relaxed in practice. For example, if $X=L^2(\D)$ and if $u \in H^2(\D)$, then $E^{\prime}(u)$ can be represented as an $L^2$-function and problem \eqref{def-Sobolev-gradient} is still well-posed. Since we will only consider $H^1$-type metrics in this work, we will not go into further detail here and refer to \cite{HenJar24} instead. To obtain the Riemannian Sobolev gradient from the Sobolev gradient, the second step involves an $X$-orthogonal projection into the tangent space $\tangentspace{u}$. For $v \in X$, the $X$-orthogonal projection $P_{u,X}: X \rightarrow \tangentspace{u}$ is given by
\begin{align}
\label{ortho_proj_def}
\left ( P_{u,X}(v) ,w\right ) _X=(v,w)_X \qquad\text{for any} \ w \in \tangentspace{u}.
\end{align}
Note that the above definition formally requires that we equip the tangent space $\tangentspace{u}$ with the $X$-metric. This is often indicated by using the notation $T_{u,X} \mathbb{S}$ instead of $T_{u} \mathbb{S}$. For the sake of readability and since the elements of $T_{u,X} \mathbb{S}$ do not change with $X$, we do not use an additional index here and just continue with $T_{u} \mathbb{S}$.
The explicit expression for the projection $P_{u,X}(v) \in \tangentspace{u}$ of $v $ is given by
\begin{align}
\label{PuX-projection}
P_{u,X}(v)=v - \frac{R_X(u)}{\|R_X(u)\|^2_X}\, (u,v)_{L^2(\mathcal{D})}.
\end{align}
Here $R_{X}(v) \in X$ is the Riesz-representation of $v$ in $X$, i.e., $ ( R_X(v),w  ) _X=(v,w)_{L^2(\mathcal{D})}$ for all $w \in \3$. 

To summarize, the Riemannian $X$-Sobolev gradient of $E$ in $u \in \mathbb{S}$ is given by  $P_{u,X} (\,\nabla_X E(u) \,)$, i.e., the $X$-orthogonal projection of the Sobolev gradient $\nabla_X E(u)$ into $\tangentspace{u}$. With these notations and operators, we can now formulate the corresponding Riemannian (conjugate) gradient method.

\subsection{The conjugate gradient method on the Riemannian manifold $\mathbb{S}$}
To introduce the concept of our method, we let $u^0 \in \mathbb{S} \subset H^1_0(\D)$ be a given starting point. For $n \geq 0$, the sequence of iterates $ \{u^{n+1} \} \subset \mathbb{S}$ of the Riemannian conjugate Sobolev gradient (RCSG) method are generated by the steps
\begin{align}
\label{riemannian-gradient-method}
u^{n+1} = \frac{ \hspace{-24pt}u^{n} + \tau_n \, d^{\, n} }{\| u^{n} + \tau_n \, d^{\, n} \|_{L^2(\D)} }.
\end{align}
That is, from $u^{n} \in \mathbb{S}$ we move by the distance $\tau_n>0$ into the search direction $d^{\, n}$. Ideally, the search direction should be a descent direction, i.e., the Riemannian gradient in direction $d^n$ is negative. The normalization $v \mapsto v / \| v \|_{L^2(\D)}$ after each iteration is the canonical retraction to $\mathbb{S}$. Unlike the Riemannian gradient method where the search direction is selected as the negative Riemannian gradient, i.e. $d^{\, n} = - P_{u^{n},X} (\,\nabla_X E(u^{n}) \,)$ (which is always a descent direction), here also the contribution from the previous search direction is included by setting
\begin{align}
\label{direction}
d^{\, n}:= \begin{cases}
- P_{u^{n},X} (\,\nabla_X E(u^{n}) \,) & \text{if } n = 0, \\
- P_{u^{n},X} (\,\nabla_X E(u^{n}) \,) + \beta^{n} P_{u^{n},X}(d^{\, n-1}) & \text{else,}
\end{cases}
\end{align}
for some momentum parameter $\beta^{n} \in \mathbb{R}$ that we specify later. Considering the fact that the Riemannian gradient $P_{u^{n},X} (\,\nabla_X E(u^{n}) \,) \in T_{u^{n}}\mathbb{S}$ and the previous search direction $d^{\, n-1} \in T_{u^{n-1}}\mathbb{S}$ both lie in different tangent spaces, we must transport the previous direction (in $T_{u^{n-1}}\mathbb{S}$) to the tangent space $T_{u^{n}}\mathbb{S}$ of the current iteration $u^n$ in order to make $d^{\, n}$ tangential to the current point $u^n$ on $\mathbb{S}$. This justifies to work with $P_{u^{n},X}(d^{\, n-1})$ instead of $d^{\, n-1}$ in \eqref{direction}. This choice is known as {\it projection-based vector transport}.

Typically, the step length $\tau_n$ in \eqref{riemannian-gradient-method} is adaptively computed to ensure optimal energy dissipation in each iteration. This can be achieved by defining:
\begin{align}
\tau_n:=\arg\min_{\tau >0}E(u^{n+1}).
\end{align}
Several optimization methods, such as Brent's method, the bisection method, the golden search section, etc., can be used to solve the above optimization problem to obtain the optimal value for $\tau_n$. It is important to note that the stability of the method is only guaranteed within a specific range of $\tau_n$.  However, providing a rigorous analytical proof of this stability requires further investigation, which is beyond the scope of this paper. For the case when uniformly $\beta^n=0$, the particular domain for searching the sequence of optimal $\tau_n$-values has been established as the interval $(0,2)$ in article \cite{PHMY24_1} for a Riemannian gradient method with adaptive metric.

As shown below, various realizations of the Riemannian conjugate gradient method can be derived from \eqref{riemannian-gradient-method} by making specific choices of (a) the inner product for the Hilbert space $X$, and (b) the real-valued momentum parameter $\beta^{n}$ in definition \eqref{direction} for the direction $d^n$. In the next subsection, we start with discussing the influence of the inner product.

\subsection{Choices for the $\boldsymbol{X}$-metric on $\mathbb{S}$}
\label{Inner_product_choices}
In the following we discuss two different inner products on the Hilbert space $\3$: the standard $H^1_0$-inner product and an energy-adaptive inner product which is defined as the linearized approximation of $E^{\prime}$ around an arbitrary linearization point $u \in \mathbb{S}$. Accordingly, our choice for the space $X$ is $(X,(\cdot,\cdot)_X)=(H^1_0(\D),(\cdot,\cdot)_X)$ equipped with the two aforementioned inner products $(\cdot,\cdot)_X$ respectively.  Since $\mathbb{S} \subset H^1_0(\D) \cong X$, this also induces a metric on the Riemannian manifold $\mathbb{S}$. For a review on more choices for the metric we refer to \cite{HenJar24}. 

\subsubsection{The $\boldsymbol{H^1_0}$-metric and the $\boldsymbol{a_u}$-metric}
We now make the choices explicit. For arbitrary $v,w \in \3$ and a linearization point $u \in \mathbb{S}$, we define
\begin{align}
\label{H1} (v,w)_{H^1_0} & :=(\nabla v, \nabla w)_{\0}, \\
\label{au} a_u(v,w) &:= (\nabla v, \nabla w)_{\0}  + (V v, w)_{\0} -\Omega (\mathcal{L}_{3} v ,w)_{\0} +\kappa(|u|^2v, w)_{\0}.
\end{align}
Note here that $a_u(u,w)=\langle E^{\prime}(u) , w \rangle$ for all $u,w \in \3$, hence we can indeed interpret $a_u(\cdot,\cdot)$ as a linearization of $E^{\prime}$. To convince ourselves that $a_u(\cdot,\cdot)$ is an inner product on $H^1_0(\D)$, we need to verify that it is positive definite. This crucially requires assumption \ref{A1} and a corresponding proof can be found in \cite[Lem.~2.2]{pc22} (for $a_0(\cdot,\cdot)$, but the result for $a_u(\cdot,\cdot)$ is directly implied). In particular, we have equivalence of the norms induced by $a_u(\cdot,\cdot)$ and $(\cdot,\cdot)_{H^1_0}$  in the sense that there exist constants $0<c \le C_u$ (that depend on $\D$, $V$, $\Omega$, $\kappa$ and, for the upper constant, on $\| u \|_{L^4(\D)}$) such that
\begin{align*}
c \,( v ,v )_{H^1_0} \,\,\le\,\, a_u(v,v) \,\,\le\,\, C_u \, ( v ,v )_{H^1_0}
\qquad \mbox{for all } v \in H^1_0(\D).
\end{align*}
Hence, both $a_u(\cdot,\cdot)$ and $(\cdot,\cdot)_{H^1_0}$ induce $H^1$-type metrics on $\mathbb{S}$ and we make the choices $(X,(\cdot,\cdot)_X)=(H^1_0(\D),a_u(\cdot,\cdot))$ and $(X,(\cdot,\cdot)_X)=(H^1_0(\D),(\cdot,\cdot)_{H^1_0})$. Note that the admissible sets of functions in $X$ are equal in both cases and that just the metric changes. 

Next, we discuss the resulting Riemannian Sobolev gradients.
%

\subsubsection{Realizations of the Riemannian $\boldsymbol{X}$-Sobolev gradients of $\boldsymbol{E}$}
In this subsection, we briefly state the Sobolev gradients associated with the $H^1_0$- and the $a_u(\cdot,\cdot)$-metric.
\begin{itemize}
\item {\bf$\mathbf{H^1_0}$-gradient:}\\
For the Hilbert space $\3$ equipped  with standard $H^1_0$-inner product, the $H^1_0$-Sobolev gradient $\nabla_{H^1_0}E(u)$ of $E$ at $u \in \3$ is given by 
\begin{align*}
    (\nabla_{H^1_0}E(u) , v)_{H^1_0} = \langle E'(u) , v \rangle  \qquad \mbox{for all } v \in \3.
\end{align*}
Recalling $E^{\prime}(u)$ from \eqref{defEprime} and that $(u,v)_{H^1_0} =(\nabla u, \nabla v)_{\0}$, we observe that 
$$ 
 \langle E'(u) , v \rangle = (u,v)_{H^1_0} + (V u - \Omega \mathcal{L}_{z} u + \kappa |u|^2 u , v)_{\0}.
 $$
This gives us
\begin{align}
\label{gradH10}
\nabla_{H^1_0}E(u) = u+ R_{H^1_0}(  V u - \Omega \mathcal{L}_{3}u + \kappa |u|^2 u),
\end{align}
where $R_{H^1_0}: \3 \rightarrow \3$ is the Ritz-projection operator given by $(R_{H^1_0}(u) , v)_{H^1_0} = (u,v)_{\0}$ for all $v \in \3$. 

Note that the practical computation of the Sobolev gradient $\nabla_{H^1_0}E(u)$ for some given $u \in \mathbb{S}$ according to \eqref{gradH10}, requires the solving of the Laplace equation of the form: Find $q \in H^1_0(\D)$ such that
\begin{align*}
- \Delta q =  V u - \Omega \mathcal{L}_{3}u + \kappa |u|^2 u.
\end{align*}
In this case we have $q=R_{H^1_0}(  V u - \Omega \mathcal{L}_{3}u + \kappa |u|^2 u)$ and consequently $\nabla_{H^1_0}E(u) = u+q$. 

From this, we obtain the corresponding Riemannian $H^1_0$-Sobolev gradient according to \eqref{PuX-projection} as 
\begin{eqnarray*}
P_{u,H^1_0} (\,\nabla_{H^1_0} E(u) \,) &=&
\nabla_{H^1_0} E(u) - \frac{R_{H^1_0}(u)}{\|R_{H^1_0}(u)\|^2_{H^1_0}}\, (u,\nabla_{H^1_0} E(u))_{L^2(\mathcal{D})} \\
&=&u+q - \frac{R_{H^1_0}(u)}{\|R_{H^1_0}(u)\|^2_{H^1_0}}\, (u,u+q)_{L^2(\mathcal{D})}, 
\end{eqnarray*}
where $R_{H^1_0}(u) \in H^1_0(\D)$ solves $- \Delta R_{H^1_0}(u) = u$. Hence, the assembly of $P_{u,H^1_0} (\,\nabla_{H^1_0} E(u) \,) $ requires the solution of two Laplace equations.

The aspect that Riemannian gradient methods based on the $H^1_0$-gradient only require the solving of the Laplace equation (with different right hand sides) can make them very attractive if they are embedded into a software package with a very efficient solver for that equation, cf. \cite{CLLZ24-discrete}.
\item {\bf$\boldsymbol{a_u(\cdot,\cdot)}$-gradient:}\\ 
Let $u \in \3$ be the linearization point and, at the same time, the point in which we want to compute the gradient of $E$. Equipping the space $\3$ with the inner product $a_{u}(\cdot,\cdot)$, we obtain the $a_{u}$-Sobolev gradient of $E$ in $u \in \mathbb{S}$ as the unique solution $\nabla_{a_u}E(u) \in H^1_0(\D)$ to
\begin{align*}
    a_u(\nabla_{a_u}E(u) , v) = \langle E'(u) , v \rangle  \qquad \mbox{for all } v \in \3.
\end{align*}
Since $\langle E'(u) , v \rangle = a_{u}(u,v)$, we readily conclude $\nabla_{a_u}E(u) = u$, i.e., the $a_u(\cdot,\cdot)$-Sobolev gradient of $E$ in $u$ is just the identity.

We can now use again \eqref{PuX-projection} to compute the corresponding Riemannian gradient, which is purely driven by the $a_u(\cdot,\cdot)$-orthogonal projection into $\tangentspace{u}$ and we have
\begin{align*}
P_{u,a_u} (\,\nabla_{a_u} E(u) \,) 
&= P_{u,a_u} (u) = u - \frac{R_{a_u}(u)}{ a_u ( R_{a_u}(u) , R_{a_u}(u) ) }\, (u,u)_{L^2(\mathcal{D})}  \\
&=  u - \frac{R_{a_u}(u)}{ ( u , R_{a_u}(u) )_{L^2(\D)} },
\end{align*}
where used $\| u\|_{L^2(\D)}=1$ and $a_u ( R_{a_u}(u) , R_{a_u}(u) )= ( u , R_{a_u}(u) )_{L^2(\D)}$. We recall that the Ritz projection $ R_{a_u}(u) \in H^1_0(\D)$ requires to solve
\begin{align*}
a_u ( R_{a_u}(u) , v ) = (u , v )_{L^2(\D)} \qquad \mbox{for all } v \in H^1_0(\D). 
\end{align*}
With this, only one linear elliptic problem has to be solved to compute the Riemannian gradient $P_{u,a_u} (\,\nabla_{a_u} E(u) \,)$. However, the differential operator changes with $u$ and has to be typically reassembled in each iteration of the (conjugate) gradient method. 

Before we conclude this subsection, we need to specify the role of the linearization point $u$ in a corresponding Riemannian gradient method. In fact, the point $u\in \mathbb{S}$ is chosen as the approximation from the previous iteration. To make this explicit, consider the corresponding gradient method that we obtain from \eqref{riemannian-gradient-method} and \eqref{direction} for $\beta^n=0$. In this case, we select the metric in iteration $n+1$ adaptively based on $u^n$ from the previous iteration. Hence, we obtain (for $\beta^n=0$)
\begin{align}
\label{dn-for-au-metric}
d^n\,\,:=\,\,-P_{u^n,a_{u^n}} (\,\nabla_{a_{u^n}} E(u^n) \,) \,\,=\,\,  u^n - \frac{R_{a_{u^n}}(u^n)}{ ( u^n , R_{a_{u^n}}(u^n) )_{L^2(\D)} }
\end{align}
with gradient step $u^{n+1} := (u^n + \tau_n \, d^n)/ \| u^n + \tau_n \, d^n \|_{L^2(\D)}$.
\end{itemize}

\subsection{Choices for the momentum parameter $\boldsymbol{\beta^n}$}
\label{choices_for_beta}
In the previous subsection, we specified two different metrics and described how the corresponding Riemannian Sobolev gradients are constructed. Next, we will turn towards the momentum parameter. In general, the momentum  parameter is chosen in such a way that the arising methods coincide with the conventional conjugate gradient method if applied to a strongly convex quadratic functional (or in other words, if applied to solve a linear problem with a symmetric, invertible operator). The specific realizations obtained e.g. by Fletcher and Reeves \cite{fletcher1964}, Polak and Ribi\'{e}re \cite{polak1969}, Dai and Yuan \cite{dai1999}, etc. were derived by considering linearizations of the original problem and by attaining convergence results for specific $\beta$ choices if the minimizing functional admits certain properties. A comprehensive abstract generalization of the parameters to Riemannian optimization problems was given by Sato \cite{Sato_2022}, which will be also the basis for our formulas for the momentum parameters.

Before presenting them, we want to illustrate the relevance of a proper $\beta$-choice as well as revealing a certain orthogonality relation for the Riemannian gradient and the previous direction which is fulfilled in our setting. In the following, we consider the RCSG method given by \eqref{riemannian-gradient-method} and \eqref{direction} for the adaptive metric $a_{u^n}(\cdot,\cdot)$ as described in the previous subsection. The specific choice of $\beta^n$ is left open for the moment.

Before we start, note that $d^n \in \tangentspace{u^n}$ is a {\it descent direction} for $E$ in $u^n$, if the Riemannian gradient is negative in direction $d^n$, which we can express as $a_{u^n}( P_{u^n,a_{u^n}} (\,\nabla_{a_{u^n}} E(u^n) \,)  , d^n )<0$. However, in the $a_{u^n}$-metric this expression simplifies tremendously by exploiting $\nabla_{a_{u^n}} E(u^n) = u^n$ and the definition of the $a_{u^n}(\cdot,\cdot)$-orthogonal projection $P_{u^n,a_{u^n}}$. Using this, we obtain that $d^n$ is a descent direction if it holds
\begin{align}
\label{descent-direction-property}
a_{u^n}( u^n , d^n )<0. 
\end{align}
With this, we return to the role of $\beta^n$. Let $u^n \in \mathbb{S}$ denote a current iterate and $d^n$ the corresponding search direction given by \eqref{dn-for-au-metric}. As mentioned before, the step length $\tau$ to reach the next iterate $u^{n+1}=\tfrac{\hspace{-19pt}u^{n}+\tau d^{n} }{ \|  u^{n}+\tau d^{n} \|  _{L^2(\mathcal{D})}}$ is determined in such a way that the function 
$$
\tau \mapsto E(u^{n+1}) := E( \tfrac{\hspace{-19pt}u^{n}+\tau d^{n} }{ \|  u^{n}+\tau d^{n} \|  _{L^2(\mathcal{D})}} )
$$ 
is minimized. To achieve this, we compute the derivative of $ E(u^{n+1})$ with respect to $\tau$ and set it to zero to determine the optimal value, i.e., we seek $\tau$ such that $\tfrac{\mathrm{d}}{\mathrm{d} \tau} E(u^{n+1}) = 0$. By the chain rule, we therefore have for the optimal $\tau$
\begin{align}
\left \langle E'(u^{n+1}), \frac{\mathrm{d}}{\mathrm{d}\tau}(u^{n+1}) \right \rangle  = 0.
\label{chainrule}
\end{align}
Noting that
\begin{align}
\frac{\mathrm{d}}{\mathrm{d} \tau} \left \|  u^{n} +\tau d^{n}\right \| _{L^2(\mathcal{D})} 
&= \frac{(u^n+\tau d^n,d^n)_{L^2(\mathcal{D})}}{\left \| u^n+\tau d^n \right \|_{L^2(\mathcal{D})} } 
\label{un_derivative}
\end{align}
we compute the derivative of $u^{n+1}$ w.r.t. $\tau$ as
\begin{eqnarray}
\nonumber \lefteqn{ \frac{\mathrm{d}}{\mathrm{d} \tau}(u^{n+1})
\,\,=\,\, \frac{\mathrm{d}}{\mathrm{d} \tau} \left(\frac{\hspace{-18pt}u^{n}+\tau d^{n}}{\left \| u^{n}+\tau d^{n}\right \| _{L^2(\mathcal{D})}}\right) 
\,\,=\,\, \frac{d^{n}\left \| u^{n} +\tau d^{n}\right \| _{L^2(\mathcal{D})} - (u^{n}+\tau d^{n})\frac{(u^n+\tau d^n,d^n)_{L^2(\mathcal{D})}}{\left \| u^n+\tau d^n \right \|_{L^2(\mathcal{D})} }}{\left \| u^{n} +\tau d^{n}\right \|^2 _{L^2(\mathcal{D})}} } \\
\nonumber  
 &=&\frac{1}{\left \| u^{n}+\tau d^{n}\right \| _{L^2(\mathcal{D})}}\left(d^n-\frac{u^n+\tau d^n}{\left \| u^{n}+\tau d^{n}\right \| _{L^2(\mathcal{D})}} \bigg(  \frac{u^n+\tau d^n}{\left \| u^{n}+\tau d^{n}\right \| _{L^2(\mathcal{D})}},d^n \bigg)_{L^2(\mathcal{D})}\right)  \hspace{70pt} \\
 \label{un1_derivative} &=& \frac{d^n - (d^n,u^{n+1})_{L^2(\mathcal{D})}u^{n+1}}{\left \| u^{n}+\tau d^{n}\right \| _{L^2}}.
\end{eqnarray}
By \eqref{un1_derivative} we can write the condition \eqref{chainrule} for the optimal step length as
\begin{align}
\left \langle E'(u^{n+1}), d^n-(d^n,u^{n+1})_{L^2(\mathcal{D})}u^{n+1} \right \rangle = 0.
\end{align}
By using the definition of the $a_u$-Sobolev gradient and the $L^2$-orthogonal projection onto the tangent space $\tangentspace{u^{n+1}}$ which is given by $P_{u^{n+1},L^2(\mathcal{D})}(d^n) = d^n-(d^n,u^{n+1})_{L^2(\mathcal{D})}u^{n+1}$, we can rewrite the condition as
\begin{align*}
a_{u^{n+1}}( \,\nabla _{a_{u^{n+1} } }  E( u^{n+1}),P_{u^{n+1},L^2(\mathcal{D}) } (d^n)  \,) =0.
\end{align*}
Since $P_{u^{n+1},L^2(\mathcal{D}) } (d^n) \in \tangentspace{u^{n+1}}$ we can exploit the $a_{u^{n+1}}(\cdot,\cdot)$-orthogonal projection onto $ \tangentspace{u^{n+1}}$ to conclude 
\begin{align}
\label{orthogonality-relation-riemannian-grad}
 a_{u^{n+1}}( \, P_{u^{n+1},a_{u^{n+1}}}(\nabla _{a_{u^{n+1}}}E(u^{n+1})) ,P_{u^{n+1},L^2}(d^n) \, ) =0.
\end{align}
Note here that $ P_{u^{n+1},a_{u^{n+1}}}(\nabla _{a_{u^{n+1}}}E(u^{n+1}))$ is nothing but the Riemannian gradient of $E$ in $u^{n+1}$ with respect to the $a_{u^{n+1}}(\cdot,\cdot)$-metric. Hence, property \eqref{orthogonality-relation-riemannian-grad} is a natural orthogonality relation between the current Riemannian gradient and the $L^2$-projection based vector transport of the previous search direction, which is fulfilled for the optimal $\tau$ value.

Thanks to our specific choice of the metric, we recall that $\nabla _{a_{u^{n+1}}}E(u^{n+1})=u^{n+1}$. Exploiting this together with $P_{u^{n+1},L^2(\mathcal{D})}(d^n)=d^n-(d^n,u^{n+1})_{L^2(\D)}u^{n+1}$ in \eqref{orthogonality-relation-riemannian-grad} and relabelling the index $n \rightarrow n-1$, we obtain
\begin{align}
\nonumber 0 &= a_{u^{n}}(P_{u^{n},a_{u^{n}}}(u^{n}),d^{n-1}-(d^{n-1},u^{n})_{L^2(\mathcal{D})}u^{n}  ) \\
\nonumber &=a_{u^{n}}( P_{u^{n},a_{u^{n}}}(u^{n}),d^{n-1}  ) -( d^{n-1}, u^{n} )_{L^2(\mathcal{D})} \,\, a_{u^{n}}(P_{u^{n},a_{u^{n}}}(u^{n}),u^{n}) \\
\label{ineq1} &\overset{\eqref{ortho_proj_def}} = a_{u^{n}}( P_{u^{n},a_{u^{n}}}(u^{n}),d^{n-1}  )-( d^{n-1}, u^{n} )_{L^2(\mathcal{D})} \, \|P_{u^{n},a_{u^{n}}}(u^{n})\|^2_{a_{u^{n}}} 
\end{align}
where for brevity $\| v \|^2_{a_{u^{n}}} := a_{u^{n}}(v,v)$. We can now identify the requirement for $\beta^n$ such that $d^n$ is a descent direction in $u^n$. For this, the Riemannian gradient $P_{u^n,a_{u^n}}(u^n)$ in the direction $d^n$ needs to be negative, i.e., according to \eqref{descent-direction-property}, we require $a_{u^n} ( u^n,d^n )<0$. To check this property, we use $d^{\, n} = - P_{u^{n},a_{u^n}} (u^{n}) + \beta^{n} P_{u^{n},a_{u^n}}(d^{\, n-1})$ to get
\begin{eqnarray}
\nonumber a_{u^n} ( u^n ,d^n )&=&- a_{u^n}( u^n, P_{u^n,a_{u^n}}(u^n) \,) +\beta^n  \, a_{u^n}( u^n, P_{u^n,a_{u^n}}(d^{n-1}) \, ) \\
\label{motivation-dai-yuan}&\overset{\eqref{ortho_proj_def}}{=}& -\left \|  P_{u^n,a_{u^n}}(u^n)\right \| ^2_{a_{u^n}}+\beta^n \,  a_{u^n}( P_{u^n,a_{u^n}} (u^n),d^{n-1} \,) \\
\nonumber &\overset{\eqref{ineq1}}{=}& \left ( -1+\beta^n(d^{n-1},u^n)_{L^2(\mathcal{D})} \right ) \left \| P_{u^{n},a_{u^{n}}}(u^{n}) \right \|^2 _{a_{u^{n}}}.
\end{eqnarray}
Consequently, $d^n$ can only be descent direction in $u^n$ if $\beta^n$ fulfills
\begin{align*}
1-\beta^n(d^{n-1},u^n)_{L^2(\mathcal{D})} > 0.
\end{align*}
The validity of this condition is crucially determined by $\beta^n$ (and implicitly through optimality by $\tau_{n-1}$). 
 A proof that shows all the momentum parameters considered here satisfy this condition is still missing in the literature and is, of course, beyond the scope of this paper. However, for the Dai-Yuan version of the parameter, the validity of the property can be made explicit, as we will demonstrate with a simple calculation below. Nevertheless, during our numerical experiments, we observed that all parameter choices produced a descent direction.

We now state the four different choices for the momentum parameter in our setting. We only formulate the parameters in the $a_{u^n}$-metric. The changes in the $H^1_0$-metric are straightforward and can be easily extracted from \cite{Sato_2022}.
%

\subsubsection{Dai--Yuan-type parameter}

Picking up on the previous discussion we can now transfer the idea by Dai and Yuan \cite{dai1999} to our setting to state the relevant momentum parameter in the energy-adaptive metric $a_{u^n}(\cdot,\cdot)$. For that, suppose that the previous search direction $d^{n-1} \in  T_{u^{n-1}}\mathbb{S}$ was a descent direction, i.e. (recalling \eqref{descent-direction-property})
\begin{align} 
\label{first}
0< - a_{u^{n-1}} ( u^{n-1}  , d^{n-1}  ) .
\end{align}
Adding 
$a_{u^n}(P_{u^{n},a_{u^n}}(u^{n}), d^{n-1}\,)$
to both sides, the descent property \eqref{first} is equivalent to
\begin{align} 
\label{descent-property-step0}
a_{u^n}(  P_{u^{n},a_{u^n}}(u^{n}) , d^{n-1} \,) <  a_{u^n}(  P_{u^n,a_{u^n}}(u^n) , d^{n-1}  )- a_{u^{n-1}}(  u^{n-1}  , d^{n-1}\big) .
\end{align}
As we saw in \eqref{motivation-dai-yuan}, the current search direction $d^{n}$ is now also a descent direction (that is $a_{u^{n}}( u^{n}  , d^{n} ) <0$) if $\beta^n$ is selected such that 
\begin{align} 
\label{descent-property-step1}
-\left \|P_{u^{n},a_{u^n}}(u^{n}) \right \|^{2}_{a_{u^{n}}} + \beta^{n}\,\, a_{u^n}(P_{u^{n},a_{u^n}}(u^{n}), d^{n-1}\,)  <0.
\end{align}
In fact, this can be fulfilled for a positive $\beta^n$ with the Dai--Yuan choice
\begin{align}
\label{eq:dy_adp}
\beta_{\mbox{\tiny DY}}^{n} := \max\left \{   0,\ \frac{\left\|P_{{u^n},a_{u^n}}(u^n)\right\|_{a_{u^n}}^2}{ a_{u^n}\big( P_{u^n,a_{u^n}}(u^n) ,  d^{n-1} ) -a_{u^{n-1}}( u^{n-1} , d^{n-1} \big)} \right \},
\end{align}
where $\| v \|_{a_{u^n}}^2 := a_{u^n}(v,v)$. To see this, let $\beta_{\mbox{\tiny DY}}^{n}>0$ (the statement is trivial for $\beta_{\mbox{\tiny DY}}^{n}=0$). In this case, the positivity of $\beta_{\mbox{\tiny DY}}^{n}$ together with the descent property \eqref{descent-property-step0} of the previous direction $d^{n-1}$ imply
\begin{align} 
\label{descent-property-step2}
\frac{a_{u^n}(P_{u^{n},a_{u^n}}(u^{n}),  d^{n-1}\,)}{ a_{u^n}( P_{u^n,a_{u^n}}(u^n) , d^{n-1} )- a_{u^{n-1}}( u^{n-1} , d^{n-1})   }  <1.
\end{align}
Hence, by using the formula \eqref{eq:dy_adp} for $\beta_{\mbox{\tiny DY}}^{n}$ we have
 \begin{eqnarray*}
\lefteqn{ -\left \|P_{u^{n},a_{u^n}}(u^{n}) \right \|^{2}_{a_{u^{n}}} + \beta^{n}\,\, a_{u^n}(P_{u^{n},a_{u^n}}(u^{n}),  d^{n-1}\,) } \\
&=& \left( \frac{a_{u^n}(P_{u^{n},a_{u^n}}(u^{n}),  d^{n-1}\,)}{ a_{u^n}( P_{u^n,a_{u^n}}(u^n) , d^{n-1} )- a_{u^{n-1}}( u^{n-1} , d^{n-1})   }  - 1\right) \, \left \|P_{u^{n},a_{u^n}}(u^{n}) \right \|^{2}_{a_{u^{n}}} \,\, \overset{\eqref{descent-property-step2}}{<}\,\, 0,
\end{eqnarray*}
which verifies \eqref{descent-property-step1} and therefore the desired descent property $a_{u^{n}}( u^{n}  , d^{n} ) <0$. This gives a justification for the choice \eqref{eq:dy_adp} in our setting. Note that the ``$\max $'' in the definition of the Dai--Yuan parameter is a consequence of the required positivity of $\beta_{\mbox{\tiny DY}}^{n}$ in our arguments. If $\beta_{\mbox{\tiny DY}}^{n}=0$, the iteration reduces to a standard Riemannian gradient step (in the $a_u$-metric). 
%

\subsubsection{Fletcher--Reeves-type parameter}
In our setting, the Fletcher--Reeves parameter can be directly extracted from the original work by Fletcher and Reeves \cite{fletcher1964}. The formula for constrained optimization problems can be also found in \cite{Smith_MIT_1998} and in particular for Riemannian optimization problems in \cite{Sato_2022}. For our $a_{u^n}$-metric, we have
\begin{align}
\beta_{\mbox{\tiny FR}}^{n} := \frac{\left \|  P_{{u^n},a_{u^n}}(u^n) \right \|_{a_{u^{n}}}^2}{\| P_{u^{n-1},a_{u^{n-1}}}(u^{n-1}) \|_{a_{u^{n-1}}}^2},
\label{eq:fr}
\end{align}
where $\| v \|_{a_{u^n}}^2 := a_{u^n}(v,v)$.

\subsubsection{Polak--Ribi\'ere-type parameter}
Following Polak and Ribi\'ere \cite{polak1969} and again Sato \cite{Sato_2022}, the corresponding momentum parameter for the $a_{u^n}$-metric is given by
\begin{align}
\beta_{\mbox{\tiny PR}}^{n} := \max\left \{ 0,\ \frac{ a_{u^{n}}( \, P_{u^n,a_{u^{n}}}(u^n), \, P_{u^n,a_{u^{n}}}(u^n)-P_{u^{n-1},a_{u^{n-1}}}(u^{n-1}) \, )}{\|  P_{u^{n-1},a_{u^{n-1}}}(u^{n-1})  \|_{a_{u^{n-1}}}^2}\right \}.
\label{eq:pr}
\end{align}

\subsubsection{Hestenes--Stiefel-type parameter}
Historically, the invention of the conjugate gradient method goes back to Hestenes and Stiefel \cite{hestenes1952}. Even though their work only considered quadratic minimization problems, their parameter choice can be heuristically transferred to more general Riemannian optimization problems, resulting in a mixture of the parameters by Polak--Ribi\'ere and Dai-Yuan. Following, the formulation provided in \cite{Sato_2022}, we obtain the Hestenes--Stiefel parameter in the $a_{u^n}$-metric as
\begin{align}
\beta_{\mbox{\tiny HS}}^{n} :=\max\left \{   0,\ \frac{ a_{u^n} ( P_{u^n,a_{u^{n}}}(u^n),P_{u^n,a_{u^{n}}}(u^n)-P_{u^{n-1},a_{u^{n-1}}}(u^{n-1}) \,) }{ a_{u^n}(  P_{u^n,a_{u^n}}(u^n) , d^{n-1} ) - a_{u^{n-1}}( u^{n-1} , d^{n-1} ) } \right \}.
\label{eq:hs}
\end{align}
In practice, for nonlinear unconstrained optimization problems, the Hestenes--Stiefel and Polak--Ribi\'ere methods typically show a similar performance and are often chosen over the Fletcher-Reeves and Dai-Yuan approaches (cf. \cite{GN1992, NoceWrig06}). We have observed the same behavior in our setting as well. Note that several authors have suggested various methods to improve computational efficiency or to reduce the number of iterations, such as restarting the scheme, i.e., repeatedly switching off the $\beta^n$-term after a certain number of steps. In our numerical experiments we did not incorporate such steps since they are mainly heuristic and require some tuning.

Looking at the standard choices for $\beta^n$ (without switching them in between) therefore ensures a fair comparison of the different methods. Note however, that we still enforce that the parameters cannot become negative as guaranteed through our definitions  \eqref{eq:dy_adp}, \eqref{eq:pr} and \eqref{eq:hs}. Only the Fletcher--Reeves parameter in \eqref{eq:fr} is positive by default. 

\section{Numerical Experiments}
\label{section:num_exp}
In this section we compare the Riemannian conjugate gradient methods from the previous section regarding their performance to compute ground states of the Gross-Pitaevskii energy functional. All our experiments were conducted on an Apple iMac-2021 with an Apple M1 Chip of 8 cores and 16 GB of RAM equipped with MATLAB 2024a. 
Our primary focus is to investigate the impact of the two different metrics (inner products) stated in Section \ref{Inner_product_choices} and the different choices for the beta parameter on the convergence speed of the RCSG method. We also compare the Riemannian conjugate gradient method with the standard Riemannian gradient method ($\beta^n=0$). One of our main observations is the notably fast convergence of the RCSG with adaptive metric ($a_{u^n}$-metric) and either Polak-Ribi\'ere or Hestenes--Stiefel momentum parameter, highlighting the great numerical efficiency of this choice in terms of a comparably low iteration number.

In the following, we consider a square domain $\mathcal{D}$=$[-L,L]^{2} \subset \mathbb{R}^{2}$, where the value of $L$ will be specified later for each experiment. 
To conduct our numerical experiments, we use a spatial discretization based on the first-order Lagrange finite elements, i.e., we consider the minimization problem \eqref{minimization-problem} over a corresponding $\mathbb{P}^1$-finite element space (on a quasi-uniform triangular mesh) and with an $L^2$-normalization constraint. Corresponding error estimates of optimal order were established in \cite{PHMY24}. The mesh size in all our experiments is $h=2*L*2^{-8}$. 
In each experiment, we used a {\it golden section search} to compute the optimal values for the step length $\tau_n$ in each iteration. The reference ground state was computed with the RCSG method with $a_{u^n}$-metric and Polak-Ribi\'ere parameter and a tolerance of $10^{-13}$ for the difference of two consecutive energies, i.e. $E(u^{n})-E(u^{n+1})$. The reason for this choice is that this combination gave us the lowest energy value compared to all other choices. 

To compare the performance of the Riemannian conjugate gradient method with different parameters, we used the stopping criterion $E(u^n)-E(u) < 10^{-9}$, where $E(u)$ denotes the energy of the reference ground state. Note that in all the experiments performed below, the interaction parameter $\kappa$, the angular velocity $\Omega$ and the trapping potential $V$ were all selected such that assumption \ref{A1} is satisfied. The trapping potential $V$ is always chosen as a harmonic potential
$$ V(x,y):=\frac{1}{2}(\gamma_x^2 x^2 + \gamma_y^2 y^2),$$ 
where the trapping frequencies $\gamma_x$ and $\gamma_y$ are specified in the individual experiments (together with $\Omega$, $\kappa$ and $L$).

\subsection{Experiment 1}
We chose the square domain with $L=6$ and the trapping frequencies, angular velocity and the interaction parameter as
$$ \gamma_x=2, \, \hspace{2mm} \gamma_y=1.9, \, \hspace{2mm} \Omega=1.9 \hspace{2mm} \text{ and} \, \hspace{2mm} \kappa=500. $$
\begin{figure}[h!]
    \flushleft
    \begin{minipage}[h]{0.5\textwidth}
        \centering
        \includegraphics[scale=0.25]{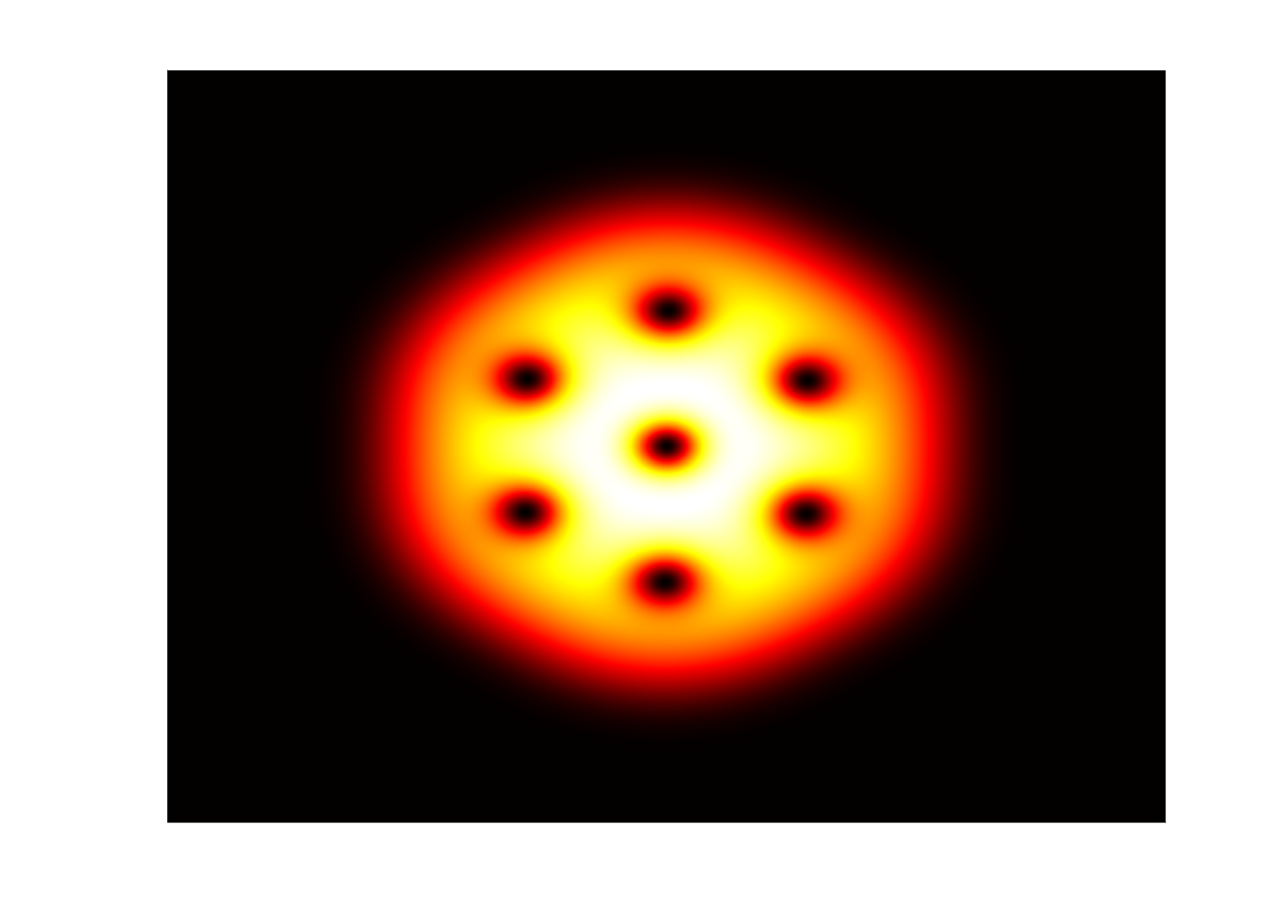}
    \end{minipage}
    \hfill 
    \begin{minipage}[h]{0.4\textwidth}
        \hspace{-13mm} 
        \centering
        \begin{tabular}{cc}
            \toprule
           \textbf{$\boldsymbol{i}$} &  \textbf{Spectrum of $\boldsymbol{E^{\prime\prime}(u)_{\vert \tangentspace{u}}}$} \\
            \midrule
            \cellcolor[gray]{.8}1 & \cellcolor[gray]{.8}20.516919568 \\
            2 & 20.517228100   \\
            3 & 20.532064567 \\
            4 & 20.533525438 \\
            5 & 20.600099691 \\
            6 & 20.614077863 \\
            7 & 20.621639553 \\
            \bottomrule
        \end{tabular}
    \end{minipage}
    \caption{{\it Left:} Surface plot of reference ground state density $|u|^2$ in 2D for experiment 1. {\it Right:} The first seven eigenvalues of $E^{\prime\prime}(u)_{\vert \tangentspace{u}}$ (ordered in ascending order based on magnitude). The ground state eigenvalue $\lambda$ appears at the bottom of the spectrum (gray).}
    \label{M2_density_plot}
\end{figure}
\begin{figure}[h!]
\centering
\begin{minipage}{0.49\textwidth}
    \centering
     \definecolor{c1}{rgb}{0.2,0.5,1.0}
     \definecolor{c2}{rgb}{0.8,0.8,0.2}
      \definecolor{c3}{rgb}{1,0.7,0.2}
    \definecolor{color0}{rgb}{0.12156862745098,0.466666666666667,0.705882352941177}
    \definecolor{color1}{rgb}{1,0.498039215686275,0.0549019607843137}
    \definecolor{color2}{rgb}{0.172549019607843,0.627450980392157,0.172549019607843}
    \definecolor{color3}{rgb}{0.83921568627451,0.152941176470588,0.156862745098039}
    \definecolor{color4}{rgb}{0.580392156862745,0.403921568627451,0.741176470588235}
    \definecolor{color5}{rgb}{0.549019607843137,0.337254901960784,0.294117647058824}
    \definecolor{color6}{rgb}{0.890196078431372,0.466666666666667,0.76078431372549}
    \begin{tikzpicture}
%
        \begin{axis}[
            name=mainaxis,
            legend style={fill=gray!20, font=\small, nodes={scale=0.6, transform shape},legend cell align=left},
            height = 01\textwidth,
            width = 1\textwidth,
            xmax   = 10000,  
            xmin   = 0, 
            ymax = 2, 
            ymin = 10^(-9),
            tick label style={font=\small}, 
            xtick={1000,3500,6000},
            scaled ticks=false,
            yticklabel style={/pgf/number format/fixed, /pgf/number format/precision=11},
            ytick={10^(-8),10^(-7),10^(-6),10^(-5),10^(-4),10^(-3),10^(-2),10^(-1)},
            xlabel=iteration $n$,
            ylabel=energy error $E(u^n)-E(u)$,    
            xlabel style={font=\small}, 
            ylabel style={font=\small},
            ymode=log
        ]
            \addplot [line width=1pt, color=c1] table {M2_DY.txt};
            \addplot[line width=1pt, color=c3] table {M2_HS.txt};
            \addplot[line width=1pt, color=color0] table {M2_FR.txt};
            \addplot[line width=1pt, color=color2] table {M2_PR.txt};
            \addplot[line width=1pt, color=color4] table {M2_RG.txt};
            \legend{
                $\text{Dai and Yuan }$,
                $\text{Hestenes and Stiefel }$,
                $\text{Fletcher and Reeves }$,
                $\text{Polak and Ribi\'ere }$,
            }
        \end{axis}
%
        \begin{axis}[
            name=insetaxis,
            at={(0.53\textwidth,0.24\textwidth)}, 
            anchor=north west,
            height = 0.40\textwidth,
            width = 0.40\textwidth,
            xmax   = 900,  
            xmin   = 400, 
            ymax = 10^(-5), 
            ymin = 10^(-9),
            tick label style={font=\tiny}, 
            xtick={600,800},
            scaled ticks=false,
            ytick={10^(-8),10^(-6)},
            ymode=log,
            clip=false
        ]
            \addplot[line width=1pt, color=c3] table {M2_HS_zoom.txt};
            \addplot[line width=1pt, color=color2] table {M2_PR_zoom.txt};
        \end{axis}
    \end{tikzpicture}
\end{minipage}
\hfill
\begin{minipage}{0.49\textwidth}
    \centering
     \definecolor{c2}{rgb}{0.2,0.8,0.6}
      \definecolor{c3}{rgb}{1,0.8,0.2}
    \definecolor{color0}{rgb}{0.12156862745098,0.466666666666667,0.705882352941177}
    \definecolor{color1}{rgb}{1,0.498039215686275,0.0549019607843137}
    \definecolor{color2}{rgb}{0.172549019607843,0.627450980392157,0.172549019607843}
    \definecolor{color3}{rgb}{0.83921568627451,0.152941176470588,0.156862745098039}
    \definecolor{color4}{rgb}{0.580392156862745,0.403921568627451,0.741176470588235}
    \definecolor{color5}{rgb}{0.549019607843137,0.337254901960784,0.294117647058824}
    \definecolor{color6}{rgb}{0.890196078431372,0.466666666666667,0.76078431372549}

    \begin{tikzpicture}
%
        \begin{axis}[
            name=mainaxis,
            legend style={fill=gray!20, font=\small, nodes={scale=0.6, transform shape},legend cell align=left},
            height = 01\textwidth,
            width = 1\textwidth,
            xmax   = 1000,  
            xmin   = 0, 
            ymax = 2, 
            ymin = 10^(-9),
            tick label style={font=\small}, 
            xtick={400,800},
            scaled ticks=false,
            yticklabel style={/pgf/number format/fixed, /pgf/number format/precision=11},
            ytick={10^(-8),10^(-7),10^(-6),10^(-5),10^(-4),10^(-3),10^(-2),10^(-1)},
            xlabel=iteration $n$,
            ylabel=energy error $E(u^n)-E(u)$,    
            xlabel style={font=\small}, 
            ylabel style={font=\small},
            ymode=log
        ]
            \addplot[line width=1pt, color=color2] table {M2_PR.txt};
            \addplot[line width=1pt, color=c2] table {M2_H1_PR.txt};
            \legend{
                $\text{Polak and Ribi\'ere with $a_{u^n}$-metric}$,
                $\text{Polak and Ribi\'ere with $H^1_0$-metric}$,
            }
        \end{axis}
 \end{tikzpicture}
\end{minipage}
    \caption{ {\it Left:} Comparison of the energy error $E(u^n)-E(u)$ per iteration for five choices of the $\beta$-parameter for Experiment 1. The manifold is adaptively equipped with the $a_{u^n}$-metric. The main plot includes an inset that shows a magnified view of the error for the Polak--Ribi\'ere and Hestenes--Stiefel parameters. {\it Right:} Comparison of the energy error $E(u^n)-E(u)$ for the Polak--Ribi\'ere parameter in the $a_{u^n}$-metric and the $H^1_0$-metric.}
    \label{M2_diff_beta}
\end{figure}
\noindent 
We initialize our method with the $L^2$-normalized interpolation of the function $u^{0}(x,y)=\tfrac{\Omega}{\sqrt{\pi}} (x+\ci y) e^{-\frac{(x^2+y^2)}{2}}$ (containing a center vortex) as suggested in \cite{BWM05}.
The  reference ground state contains quantized vortices, and the corresponding density $|u|^2$ is depicted in Figure \ref{M2_density_plot} (left). 
The reference ground state has an energy of $E(u)\approx 7.168961589$ and the corresponding ground state eigenvalue is $\lambda \approx 20.516919568$ (for minimizers in the selected finite element space). Recalling the second order optimality condition, we also verified that the ground state eigenvalue $\lm$ is the smallest (and simple) eigenvalue of $E''(u) \vert_{\tangentspace{u}}$, as shown in the table of Figure \ref{M2_density_plot}. As expected, the corresponding eigenfunction of $E''(u) \vert_{\tangentspace{u}}$ to $\lambda$ is given by $\ci u$. 
This confirms that the sufficient condition for $u$ being a (quasi-isolated) local minimizer is satisfied for our reference ground state.

To compare the performance of the RCSG methods \eqref{riemannian-gradient-method}-\eqref{direction} with adaptive inner product $(\cdot ,\cdot)_{a_{u^n}}$ and four different choices of the momentum parameter as stated in Subsection \ref{choices_for_beta}, along with the standard choice of $\beta^n=0$, we plot the corresponding numbers of iterations versus the error $E(u^n)-E(u)$ in Figure \ref{M2_diff_beta} (left).
The plot shows that the Polak--Ribi\'ere and Hestenes--Stiefel parameters perform similarly, and both perform significantly better than the other two available options, and, as expected, also better than the Riemannian gradient method ($\beta^n=0$), which took nearly 10,000 iterations to reach an energy error of order $10^{-3}$. 

It should be noted that in order to prevent degeneracy of the step size $\tau_n$, we enforced that $\tau_n$ should not fall below a value of $0.001$. If there is no $\tau_n \ge 0.001$ such that the energy reduces in the $n$'th iteration, we default to $\beta^n=0$ and hence take a step with the standard Riemannian gradient method for which $\tau_n$ cannot degenerate (cf. \cite{PHMY24_1}).
As a result of this strategy, for very few iterations with the Dai--Yuan and Fletcher--Reeves parameters, we had to switch to $\beta^n=0$ in order to avoid any increase in energy. We did not observe this effect for the Polak--Ribi\'ere and Hestenes--Stiefel parameters.

 Finally, we also compared the performance of the RCSG method with Polak--Ribi\'ere parameter in the adaptive metric with the corresponding method in the standard $H^1_0$-metric \eqref{H1}. The plot is presented in Figure \ref{M2_diff_beta} (right). Although in this particular experiment the RCSG method with (adaptive) $a_{u^n}$-metric  only required 100 iterations less than the realization with $H^1_0$-metric (which is not much considering the total number of iterations), the difference in performance became much more pronounced in our other experiments, where the methods with $H^1_0$-metric required a significantly larger number of iterations. However, we also noticed that all methods converged to the same ground state in Experiment 1, regardless of the choice of the $\beta^n$ parameter in the search direction.  We will revisit this perspective in Experiment 3, where we observed that this is not always the case.


\subsection{Experiment 2}
\begin{figure}[h!]
    \flushleft
    \begin{minipage}[h]{0.5\textwidth}
        \centering
        \includegraphics[scale=0.25]{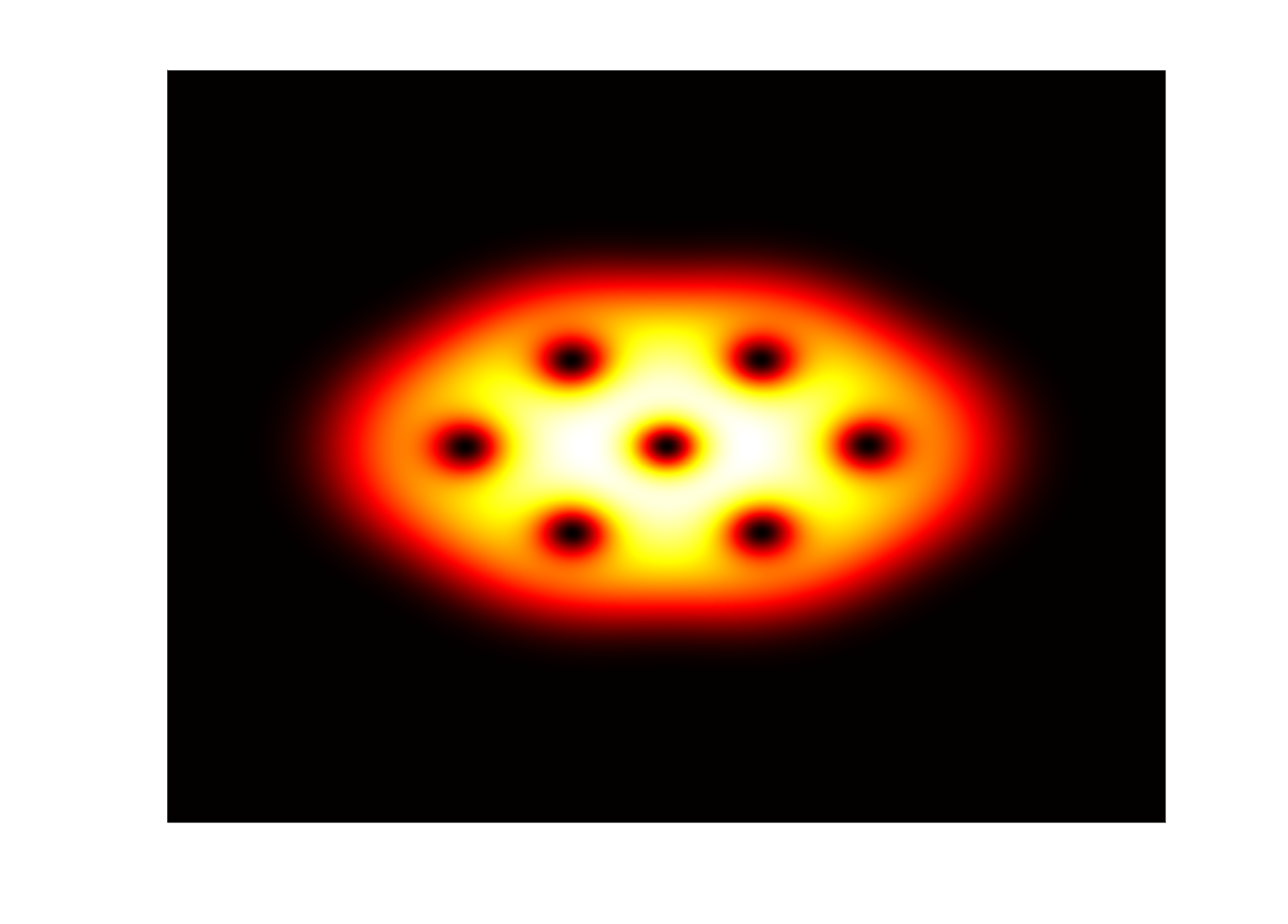}
    \end{minipage}
    \hfill 
    \begin{minipage}[h]{0.4\textwidth}
    \hspace{-13mm} 
        \centering
        \begin{tabular}{cc}
            \toprule
           \textbf{$\boldsymbol{i}$} &  \textbf{Spectrum of $\boldsymbol{E^{\prime\prime}(u)_{\vert \tangentspace{u}}}$} \\
            \midrule
            \cellcolor[gray]{.8}1 & \cellcolor[gray]{.8}10.994845147  \\
            2 & 10.996312407   \\
            3 & 10.999259973 \\
            4 & 10.999977264 \\
            5 & 11.007689478 \\
            6 & 11.017599133 \\
            7 &  11.023446235 \\
            \bottomrule
        \end{tabular}
    \end{minipage}
    \caption{{\it Left:} Surface plot of reference ground state density $|u|^2$ for Experiment 2. {\it Right:} The first seven eigenvalues of $E^{\prime\prime}(u)_{\vert \tangentspace{u}}$ (ordered in ascending order based on magnitude). The ground state eigenvalue $\lambda$ appears at the bottom of the spectrum (gray).}
    \label{M1_density_plot}
\end{figure}

We select $L=8$, i.e. $\mathcal{D}$=$[-8,8]^2$, and the parameters for trapping frequencies, angular velocity and interaction between the particles as $ \gamma_x=1.1, \, \hspace{2mm} \gamma_y=1.3, \, \hspace{2mm} \Omega=1.2 \hspace{2mm}  \text{ and} \, \hspace{2mm} \kappa=400. $
We initialize the methods with the complex conjugate of the $L^2$-normalized interpolation of the function $u^0(x,y)$ as in Experiment 1. The energy of the reference ground state is computed as $E(u)\approx 3.886043618$ and the corresponding ground state eigenvalue as $\lambda \approx 10.994845147 $. The plot of the corresponding density is shown in Figure \ref{M1_density_plot},  and the table therein shows that $u$ is at least a local minimizer (as explained in Section \ref{Model_setting} and for Experiment 1).  The plot in Figure \ref{M1_diff_beta} (left) compares the performance of all five different momentum parameters for the $a_{u^n}$-metric.  Again, the Polak--Ribi\'ere and  Hestenes--Stiefel parameters perform better than the other three possibilities by a substantial margin. Moreover, the plot in Figure \ref{M1_diff_beta} (right) shows that the method in the $a_{u^n}$-metric requires significantly less iterations than the method in the $H^1_0$-metric.
%

\begin{figure}[h!]
    \centering
    \begin{minipage}{0.48\textwidth}
        \centering
         \definecolor{c1}{rgb}{0.2,0.5,1.0}
     \definecolor{c2}{rgb}{0.8,0.8,0.2}
      \definecolor{c3}{rgb}{1,0.7,0.2}
        \definecolor{color0}{rgb}{0.12156862745098,0.466666666666667,0.705882352941177}
        \definecolor{color1}{rgb}{1,0.498039215686275,0.0549019607843137}
        \definecolor{color2}{rgb}{0.172549019607843,0.627450980392157,0.172549019607843}
        \definecolor{color3}{rgb}{0.83921568627451,0.152941176470588,0.156862745098039}
        \definecolor{color4}{rgb}{0.580392156862745,0.403921568627451,0.741176470588235}
        \definecolor{color5}{rgb}{0.549019607843137,0.337254901960784,0.294117647058824}
        \definecolor{color6}{rgb}{0.890196078431372,0.466666666666667,0.76078431372549}
        \begin{tikzpicture}
            \begin{axis}[
                name=mainaxis,
                legend style={fill=gray!20, font=\small, nodes={scale=0.6, transform shape},legend cell align=left},
                height = 01\textwidth,
                width = 1\textwidth,
                xmax   = 10000,  
                xmin   = 0, 
                ymax = 2, 
                ymin = 10^(-9),
                tick label style={font=\small}, 
                xtick={1000,5000,9000},
                scaled ticks=false,
                yticklabel style={/pgf/number format/fixed, /pgf/number format/precision=11},
                ytick={10^(-8),10^(-7),10^(-6),10^(-5),10^(-4),10^(-3),10^(-2),10^(-1)},
                xlabel=iteration $n$,
                ylabel=energy error $E(u^n)-E(u)$,    
                xlabel style={font=\small}, 
                ylabel style={font=\small},
                ymode=log
            ]
                \addplot [line width=1pt, color=c1] table {M1_YD_cut.txt};
                \addplot[line width=1pt, color=c3] table {M1_HS.txt};
                \addplot[line width=1pt, color=color0] table {M1_FR.txt};
                \addplot[line width=1pt, color=color2] table {M1_PR.txt};
                \addplot[line width=1pt, color=color4] table {M1_RG_cut.txt};
                \legend{
                    $\text{Dai and Yuan}$,
                    $\text{Hestenes and Stiefel}$,
                    $\text{Fletcher and Reeves}$,
                     $\text{Polak and Ribi\'ere }$,
                    $\text{Riemannian Gradient}$,                    
                }
            \end{axis}

            \begin{axis}[
                name=insetaxis,
                at={(0.2\textwidth,0.24\textwidth)}, 
                anchor=north west,
                height = 0.40\textwidth,
                width = 0.40\textwidth,
                xmax   = 300,  
                xmin   = 0, 
                ymax = 2, 
                ymin = 10^(-9),
                tick label style={font=\tiny}, 
                xtick={100,200},
                scaled ticks=false,
                ytick={10^(-6),10^(-3)},
                ymode=log,
                clip=false
            ]
                \addplot[line width=1pt, color=c3] table {M1_HS.txt};
                \addplot[line width=1pt, color=color2] table {M1_PR.txt};
            \end{axis}
        \end{tikzpicture}
    \end{minipage}
    \hfill
    \begin{minipage}{0.48\textwidth}
        \centering
         \definecolor{c1}{rgb}{0.2,0.5,1.0}
     \definecolor{c2}{rgb}{0.2,0.8,0.6}
      \definecolor{c3}{rgb}{1,0.7,0.2}
        \definecolor{color0}{rgb}{0.12156862745098,0.466666666666667,0.705882352941177}
        \definecolor{color1}{rgb}{1,0.498039215686275,0.0549019607843137}
        \definecolor{color2}{rgb}{0.172549019607843,0.627450980392157,0.172549019607843}
        \definecolor{color3}{rgb}{0.83921568627451,0.152941176470588,0.156862745098039}
        \definecolor{color4}{rgb}{0.580392156862745,0.403921568627451,0.741176470588235}
        \definecolor{color5}{rgb}{0.549019607843137,0.337254901960784,0.294117647058824}
        \definecolor{color6}{rgb}{0.890196078431372,0.466666666666667,0.76078431372549}
        \begin{tikzpicture}
            \begin{axis}[
                name=mainaxis,
                legend style={fill=gray!20, font=\small, nodes={scale=0.6, transform shape},legend cell align=left},
                height = 01\textwidth,
                width = 1\textwidth,
                xmax   = 2704,  
                xmin   = 0, 
                ymax = 2, 
                ymin = 10^(-9),
                tick label style={font=\small}, 
                xtick={500,1500,2500},
                scaled ticks=false,
                yticklabel style={/pgf/number format/fixed, /pgf/number format/precision=11},
                ytick={10^(-8),10^(-7),10^(-6),10^(-5),10^(-4),10^(-3),10^(-2),10^(-1)},
                xlabel=iteration $n$,
                ylabel=energy error $E(u^n)-E(u)$,    
                xlabel style={font=\small}, 
                ylabel style={font=\small},
                ymode=log
            ]
                \addplot[line width=1pt, color=color2] table {M1_PR.txt};
                \addplot[line width=1pt, color=c2] table {M1_H1_PR.txt};
                \legend{
                    $\text{Polak and Ribi\'ere with $a_{u^n}$-metric}$,
                    $\text{Polak and Ribi\'ere with }H^1$-$\text{metric}$,
                }
            \end{axis}
        \end{tikzpicture}
    \end{minipage}
     \caption{ Comparision of the energy error $E(u^n)-E(u)$ within the first $10^4$ iterations for Experiment 2. {\it Left}: In all methods, the manifold is adaptively equipped with the $a_{u^n}$-metric. The main plot includes an inset that shows a magnified view of the error for the Polak--Ribi\'ere and Hestenes--Stiefel parameters. {\it Right}: Comparison of the RCSG methods with Polak--Ribi\'ere parameter in the $a_{u^n}$-metric and the $H^1_0$-metric.}
        \label{M1_diff_beta}
\end{figure}


\subsection{Experiment 3: Stagnation in neighborhoods of excited states}
\begin{figure}[h!]
    \flushleft
    \begin{minipage}[h]{0.5\textwidth}
        \centering
        \includegraphics[scale=0.25]{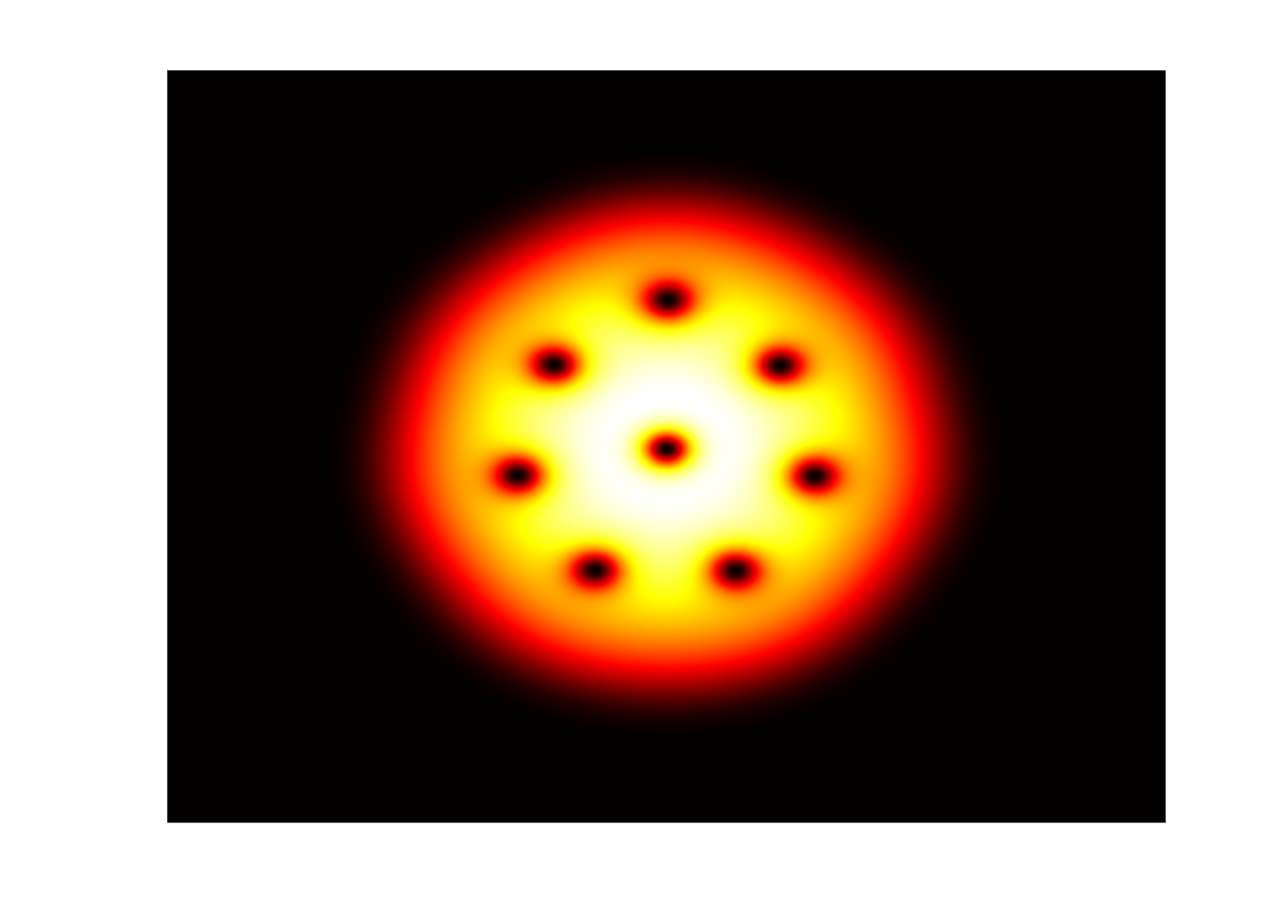}
    \end{minipage}
    \hfill 
    \begin{minipage}[h]{0.4\textwidth}
    \hspace{-13mm} 
        \centering
        \begin{tabular}{cc}
            \toprule
           \textbf{$\boldsymbol{i}$} &  \textbf{Spectrum of $\boldsymbol{E^{\prime\prime}(u)_{\vert \tangentspace{u}}}$} \\
            \midrule
            \cellcolor[gray]{.8}1 & \cellcolor[gray]{.8}162.8748496  \\
            2 & 162.8756237   \\
            3 & 162.9536145 \\
            4 & 162.9561402 \\
            5 & 163.0294160 \\
            6 & 163.1439511 \\
            7 & 163.2677753 \\
            \bottomrule
        \end{tabular}
    \end{minipage}
    \caption{{\it Left:} Surface plot of reference ground state density $|u|^2$ for Experiment 3. {\it Right:} The first seven eigenvalues of $E^{\prime\prime}(u)_{\vert \tangentspace{u}}$ (ordered in ascending order based on magnitude). The ground state eigenvalue $\lambda$ appears at the bottom of the spectrum (gray).}
    \label{M3_density_plot}
\end{figure}

We let $L=3$ and select the data parameters as $ \gamma_x=11, \, \hspace{1mm} \gamma_y=10, \, \hspace{1mm} \Omega=9 \, \text{ and} \, \hspace{1mm} \kappa=1000$.
As suggested in \cite{BWM05}, the methods are initialized with an interpolation of the following $L^2$-normalized function, 
$$ 
u^{0}(x,y)=\tfrac{\Omega}{\sqrt{\pi}} (x+\ci y) e^{-\frac{(x^2+y^2)}{2}} + \tfrac{1-\Omega}{\sqrt{\pi}}e^{-\frac{(x^2+y^2)}{2}}.$$ 
The energy of the reference ground state (in the selected finite element space) is determined as $E\approx 55.92566166$ and the corresponding ground state eigenvalue as $\lm \approx 162.8748496$. The presence of the ground state eigenvalue at the bottom of the spectrum of $E^{\prime\prime}(u)_{\vert \tangentspace{u}}$ proves that the computed $u$ is a local minimizer (see Figure \ref{M3_density_plot}). Again, in comparing the performance of our method in the $a_{u^n}$-metric for different momentum parameters, the Polak--Ribi\'ere and Hestenes--Stiefel parameters outperformed all the others, see Figure \ref{M3_diff_beta} (left).
Moreover, in Figure \ref{M3_diff_beta} (right) we also observe a significant reduction of the iteration numbers when the Riemannian conjugate gradient method is used in combination with the $a_{u^n}$-metric compared to the $H^1_0$-metric.
%

\begin{figure}[h!]
    \centering
    \begin{minipage}{0.48\textwidth}
        \centering
         \definecolor{c1}{rgb}{0.1,0.6,0.9}
     \definecolor{c2}{rgb}{0.8,0.8,0.2}
      \definecolor{c3}{rgb}{1,0.7,0.2}
        \definecolor{color0}{rgb}{0.12156862745098,0.466666666666667,0.705882352941177}
        \definecolor{color1}{rgb}{1,0.498039215686275,0.0549019607843137}
        \definecolor{color2}{rgb}{0.172549019607843,0.627450980392157,0.172549019607843}
        \definecolor{color3}{rgb}{0.83921568627451,0.152941176470588,0.156862745098039}
        \definecolor{color4}{rgb}{0.580392156862745,0.403921568627451,0.741176470588235}
        \definecolor{color5}{rgb}{0.549019607843137,0.337254901960784,0.294117647058824}
        \definecolor{color6}{rgb}{0.890196078431372,0.466666666666667,0.76078431372549}

        \begin{tikzpicture}
            \begin{axis}[
                name=mainaxis,
                legend style={fill=gray!20, font=\small, nodes={scale=0.6, transform shape},  legend cell align=left, at={(0.95,0.53)},},
                height = 01\textwidth,
                width = 01\textwidth,
                xmax   = 2198,  
                xmin   = 0, 
                ymax = 30, 
                ymin = 10^(-9),
                tick label style={font=\small}, 
                xtick={500,1500},
                scaled ticks=false,
                yticklabel style={/pgf/number format/fixed, /pgf/number format/precision=11},
                ytick={10^(-8),10^(-7),10^(-6),10^(-5),10^(-4),10^(-3),10^(-2),10^(-1)},
                xlabel=iteration $n$,
                ylabel=energy error $E(u^n)-E(u)$,    
                xlabel style={font=\small}, 
                ylabel style={font=\small},
                ymode=log
            ]
             \addplot[line width=1pt, color=c1] table {M3_DY_2198.txt};
                \addplot[line width=1pt, color=c3] table {M3_HS.txt};
                \addplot[line width=1pt, color=color0] table {M3_FR.txt};
                \addplot[line width=1pt, color=color2] table {M3_PR.txt}; 
                \addplot[line width=1pt, color=color4] table {M3_RG_2198.txt};
                \legend{
                   $\text{Dai and Yuan }$,
                   $\text{Hestenes and Stiefel }$,
                   $\text{Fletcher and Reeves }$,
                   $\text{Polak and Ribi\'ere }$,
                   $\text{Riemannian Gradient }$,
                }
            \end{axis}

            \begin{axis}[
                name=insetaxis,
                at={(0.45\textwidth,0.07\textwidth)}, 
                height = 0.35\textwidth,
                width = 0.45\textwidth,
                xmax   = 546,  
                xmin   = 542, 
                ymax = 2*10^(-9), 
                ymin = 8.5*10^(-10),
                tick label style={font=\tiny}, 
                xtick={544},
                scaled ticks=false,
                ymode=log,
                clip=false
            ]
                \addplot[line width=1pt, color=c3] table {M3_HS_zoom.txt};
                \addplot[line width=1pt, color=color2] table {M3_PR_zoom.txt};
            \end{axis}
        \end{tikzpicture}
    \end{minipage}
    \hfill
    \begin{minipage}{0.48\textwidth}
        \centering
         \definecolor{c1}{rgb}{0.2,0.5,1.0}
    \definecolor{c2}{rgb}{0.2,0.8,0.6}
      \definecolor{c3}{rgb}{1,0.7,0.2}
        \definecolor{color0}{rgb}{0.12156862745098,0.466666666666667,0.705882352941177}
        \definecolor{color1}{rgb}{1,0.498039215686275,0.0549019607843137}
        \definecolor{color2}{rgb}{0.172549019607843,0.627450980392157,0.172549019607843}
        \definecolor{color3}{rgb}{0.83921568627451,0.152941176470588,0.156862745098039}
        \definecolor{color4}{rgb}{0.580392156862745,0.403921568627451,0.741176470588235}
        \definecolor{color5}{rgb}{0.549019607843137,0.337254901960784,0.294117647058824}
        \definecolor{color6}{rgb}{0.890196078431372,0.466666666666667,0.76078431372549}

        \begin{tikzpicture}
            \begin{axis}[
                name=mainaxis,
                legend style={fill=gray!20, font=\small, nodes={scale=0.6, transform shape},legend cell align=left},
                height = 01\textwidth,
                width = 01\textwidth,
                xmax   = 10000,  
                xmin   = 0, 
                ymax = 2, 
                ymin = 10^(-9),
                tick label style={font=\small}, 
                xtick={1000,5000,9000},
                scaled ticks=false,
                yticklabel style={/pgf/number format/fixed, /pgf/number format/precision=11},
                ytick={10^(-8),10^(-7),10^(-6),10^(-5),10^(-4),10^(-3),10^(-2),10^(-1)},
                xlabel=iteration $n$,
                ylabel=energy error $E(u^n)-E(u)$,    
                xlabel style={font=\small}, 
                ylabel style={font=\small},
                ymode=log
            ]
                \addplot[line width=1pt, color=color2] table {M3_PR.txt};
                \addplot[line width=1pt, color=c2] table {M3_H1_PR.txt};
                \legend{
                   $\text{Polak and Ribi\'ere with $a_{u^n}$-metric}$,
                   $\text{Polak and Ribi\'ere with H1-metric}$,
                }
            \end{axis}
        \end{tikzpicture}
         \caption{Comparison of the energy error $E(u^n)-E(u)$ for experiment 3. {\it Left}: The manifold is adaptively equipped with the $a_{u^n}$-metric. The main plot includes an inset that shows a magnified view of the error for the Polak--Ribi\'ere and Hestenes--Stiefel parameter. {\it Right}: Comparison for the Polak--Ribi\'ere parameter in the $a_{u^n}$-metric and the $H^1_0$-metric.}
    \label{M3_diff_beta}
    \end{minipage}
\end{figure}

\begin{figure}[h!]
    \begin{minipage}[h]{0.45\textwidth}
        \centering
        \includegraphics[scale=0.25]{ground_state_exp3-eps-converted-to.pdf} 
         \centering
        \hspace{8mm}$E=55.925661, \, \lm=162.87484$
    \end{minipage}
    \hfill 
    \hspace{-5mm}
     \begin{minipage}[h]{0.45\textwidth}
        \centering
        \includegraphics[scale=0.25]{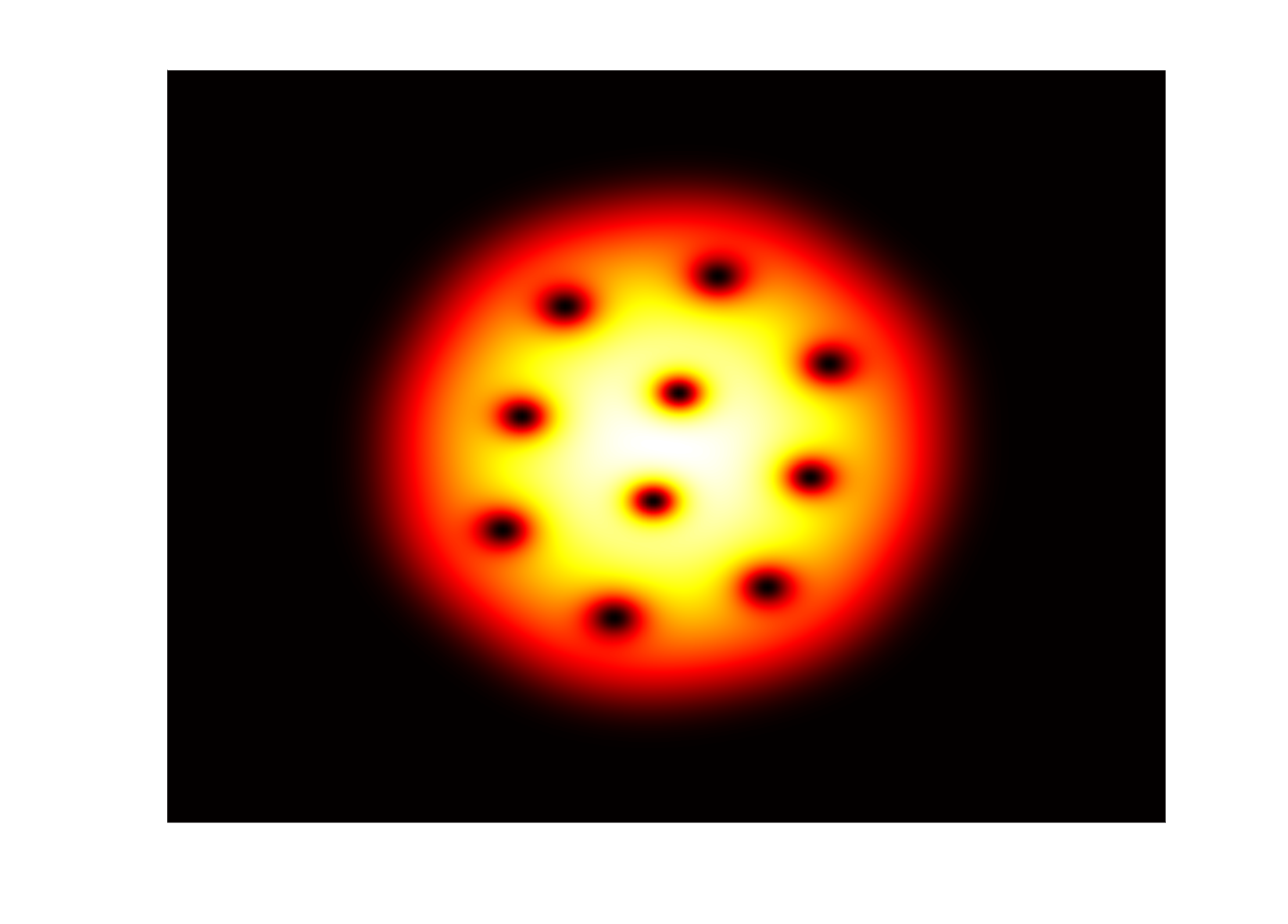}
         \centering
         $E=55.987960, \, \lm=161.92972$
     \end{minipage}
    \caption{Surface plots of the density functions $|u|^2$ for Experiment 3. {\it Left:} Density of the reference ground state. {\it Right:} Density of the computed excited state (saddle point of $E$ on $\mathbb{S}$).}
    \label{M3_ES}
\end{figure}

We also observed that the Riemannian gradient method ($\beta^n=0$) did not convergence to the ground state, but got stuck near an excited state. This did not only happen for Experiment 3, but was also noticed for other settings not included in this article. The given experiment hence serves as an illustrative example to observe this drawback. To be more specific, we noticed that even after more than 10,000 iterations of the Riemannian gradient method (for $\beta^n=0$), the difference between the energy of the converged state $u_{\mbox{\tiny e}}$ and the reference ground state $u$ is on the order of $10^{-1}$. We also computed the corresponding residual $E^{\prime}(u_{\mbox{\tiny e}})- \lambda_{\mbox{\tiny e}} I u_{\mbox{\tiny e}}$ and the spectrum of $E^{\prime\prime}(u_{\mbox{\tiny e}})\vert_{\tangentspace{u_{\mbox{\tiny e}}}}$ to confirm that the converged state $u_{\mbox{\tiny e}}$ is indeed a saddle point of $E$ on $\mathbb{S}$ and not a local minimizer.
The corresponding densities of the reference ground state $u$ and the excited state $u_{\mbox{\tiny e}}$ are plotted in Figure \ref{M3_ES}. An important point to notice here is that the eigenvalue corresponding to the excited state found by the Riemannian gradient method is smaller than the ground state eigenvalue. Hence, the ground state eigenvalue is not necessarily the smallest eigenvalue of the Gross-Pitaevskii eigenvalue problem \eqref{eigen_value_problem_1}.

%
\subsection{Experiment 4: Number of iterations versus the angular frequency $\Omega$}
In the last experiment, we examine how the number of iterations is affected by a growing number of vortices. More vortices can be achieved by increasing the angular velocity $\Omega$ (which has to stay below a critical limit) while keeping the other parameters fixed. Here we select $\mathcal{D} = [-6,6]^2$, the trapping potential is chosen as $V(x,y) = 2(x^2 + y^2)$, and the interaction parameter $\kappa = 200$. We study $\Omega=2.0,\,\,2.4,\,\,2.7$. To ensure that the finite element space is rich enough to resolve the vortices, we deviate from the previous experiments and chose a spatial discretization based on $\mathbb{P}^2$ finite elements on a simplicial mesh of mesh width $h = 12 * 2^{-8}$. We calculate the reference ground state with the same stopping criterion as in earlier experiments along with the initial condition as an $L^2$-normalized function:
$$ 
u^{0}(x,y)=\tfrac{\Omega}{\sqrt{\pi}} (x+\ci y) e^{-\frac{(x^2+y^2)}{2}} + \tfrac{1-\Omega}{\sqrt{\pi}}e^{-\frac{(x^2+y^2)}{2}}.$$ 
After the stopping criterion is met, the energy of the ground state (at least local minimizer) in each setting is observed as follows: for $\Omega =2$, $E(u)=4.96118965$; for $\Omega = 2.4$, $E(u)=4.19072815$; and for $\Omega=2.7$, $E(u)=3.04745845$. To investigate how increasing $\Omega$ values affect the iteration numbers, we restrict the experiment to the Polak--Ribi\'ere-realization of the Riemannian conjugate gradient method as it showed the best performance in the previous experiments (together with the Hestenes--Stiefel realization). 
Figure \ref{diff_omega} shows the corresponding iteration numbers versus the error $E(u^n)-E(u)$ for the different $\Omega$ values. As evident from the figure, the number of iterations increases significantly with higher angular velocities. This demonstrates that for fast-rotating condensates, which contain a greater number of vortices, the computational complexity increases considerably compared to slowly rotating condensates. The increased complexity likely arises due to the intricate vortex structures and the associated energy landscape, which requires more iterations for the algorithm to converge. In these cases, the usage of an efficient metric-driven gradient method becomes even more important.

\begin{figure}[h!]
    \begin{minipage}[h]{0.30\textwidth}
        \centering
        \includegraphics[scale=0.2]{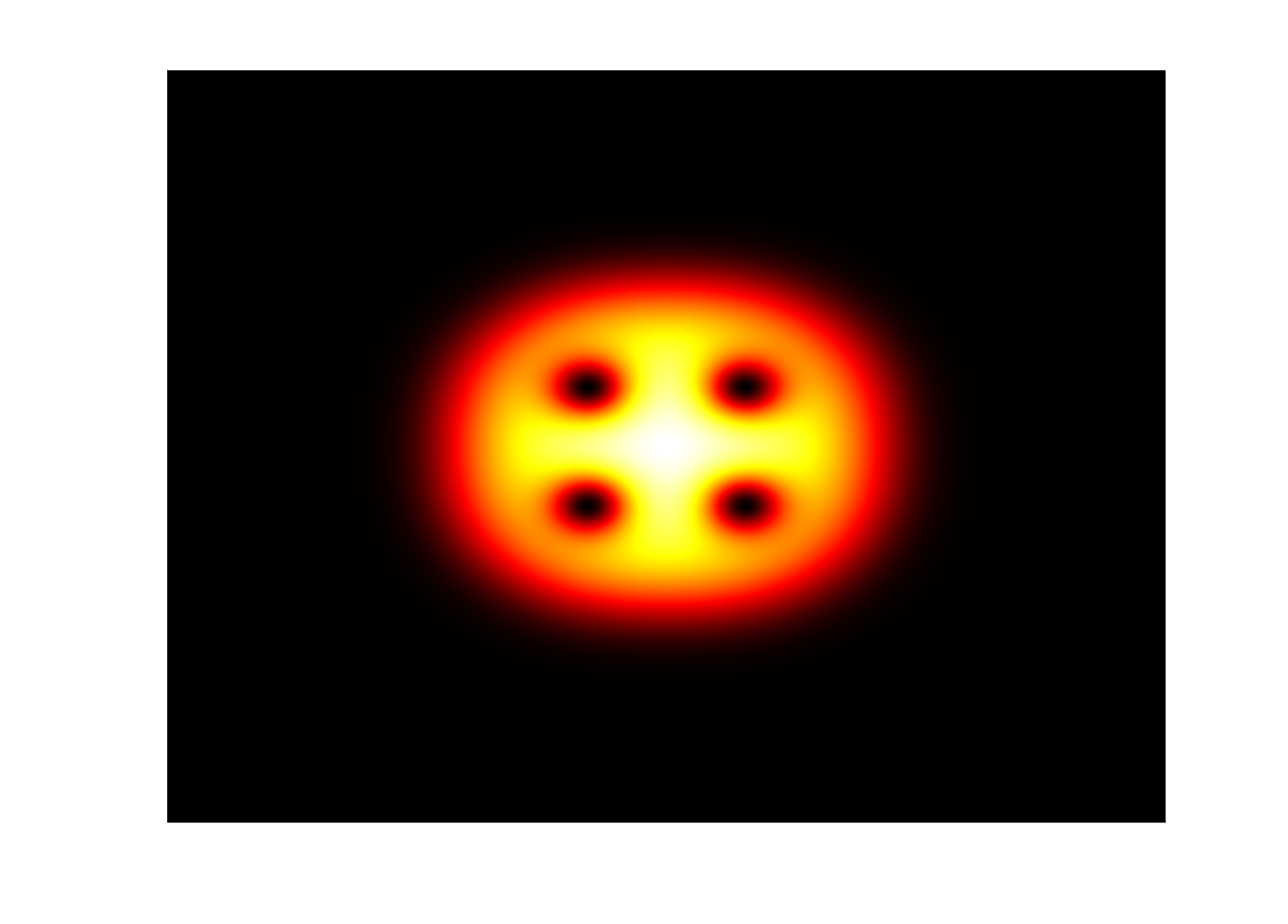} 
         \centering
    \end{minipage}
    \hfill 
    \hspace{-5mm}
     \begin{minipage}[h]{0.30\textwidth}
        \centering
        \includegraphics[scale=0.2]{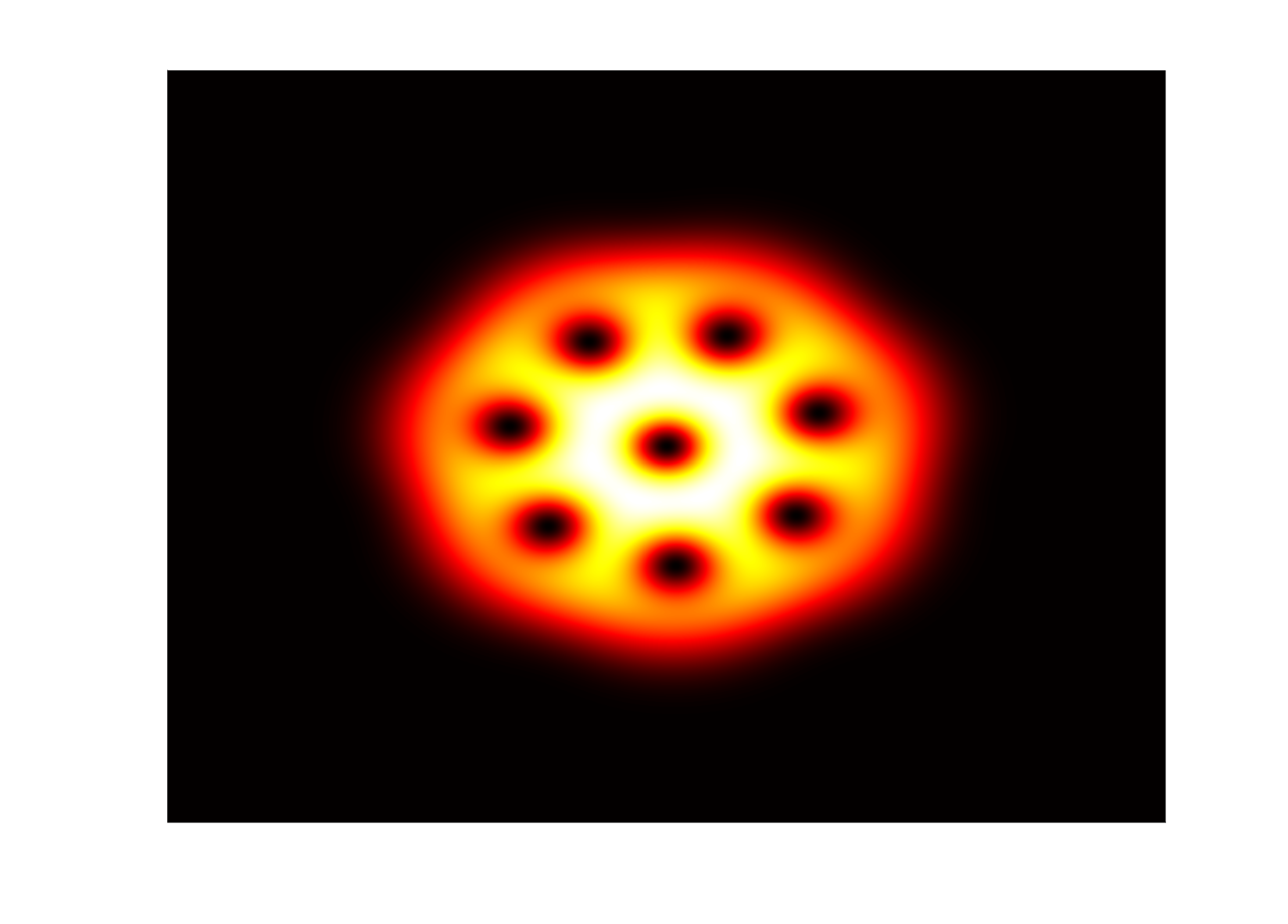}
         \centering
     \end{minipage}
      \hfill 
    \hspace{-5mm}
     \begin{minipage}[h]{0.30\textwidth}
        \centering
        \includegraphics[scale=0.2]{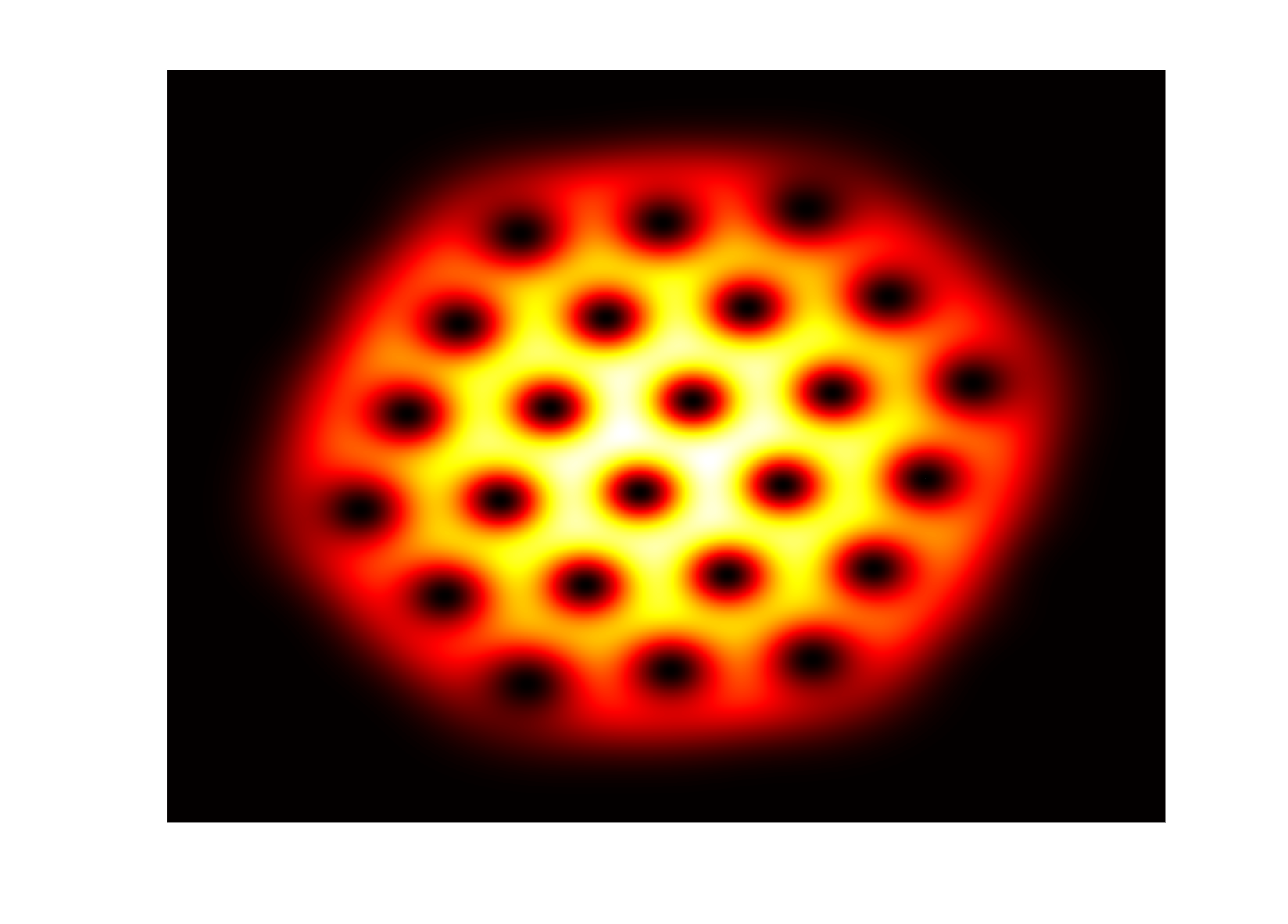} 
         \centering
     \end{minipage}
    \caption{Surface plots of the density functions $|u|^2$ for Experiment 4 with $\Omega=2$ (left), $\Omega=2.4$ (middle) and $\Omega=2.7$ (right).}
    \label{exp4_density}
\end{figure}
%
%


\begin{figure}[h!]
    \centering
    \hfill
         \definecolor{c1}{rgb}{0.2,0.5,1.0}
    \definecolor{c2}{rgb}{0.2,0.8,0.6}
      \definecolor{c3}{rgb}{1,0.7,0.2}
        \definecolor{color0}{rgb}{0.12156862745098,0.466666666666667,0.705882352941177}
        \definecolor{color1}{rgb}{1,0.498039215686275,0.0549019607843137}
        \definecolor{color2}{rgb}{0.172549019607843,0.627450980392157,0.172549019607843}
        \definecolor{color3}{rgb}{0.83921568627451,0.152941176470588,0.156862745098039}
        \definecolor{color4}{rgb}{0.580392156862745,0.403921568627451,0.741176470588235}
        \definecolor{color5}{rgb}{0.549019607843137,0.337254901960784,0.294117647058824}
        \definecolor{color6}{rgb}{0.890196078431372,0.466666666666667,0.76078431372549}
\begin{tikzpicture}
            \begin{axis}[
                name=mainaxis,
                legend style={fill=gray!20, font=\small, nodes={scale=0.6, transform shape},legend cell align=left},
                height = 0.4\textwidth,
                width = 01\textwidth,
                xmax   = 1300,  
                xmin   = 0, 
                ymax = 2, 
                ymin = 10^(-9),
                tick label style={font=\small}, 
                xtick={126,467,938},
                scaled ticks=false,
                yticklabel style={/pgf/number format/fixed, /pgf/number format/precision=11},
                ytick={10^(-8),10^(-6),10^(-4),10^(-2)},
                xlabel=iteration $n$,
                ylabel= energy error $E(u^n)-E(u)$,    
                xlabel style={font=\small}, 
                ylabel style={font=\small},
                ymode=log
            ]
                \addplot[line width=1.2pt, color=c2] table {Compare2_check1.txt};
                \addplot[line width=1.2pt, color=c3] table {Compare24_check1.txt};
                  \addplot[line width=1.2pt, color=color4] table {Compare27_check1.txt};
                \legend{
                   $\Omega = 2.0$,
                   $\Omega = 2.4$,
                    $\Omega = 2.7$,
                }
            \end{axis}
        \end{tikzpicture}
         \caption{Comparison of the energy error $E(u^n)-E(u)$ in Experiment 4 obtained with the Polak--Ribi\'ere RCSG method for different angular velocities $\Omega$.}         
    \label{diff_omega}
   \end{figure}

\section{Conclusions}
In this paper we proposed a class of metric-driven Riemannian conjugate gradient methods which are based on the concept of Sobolev gradients with an energy-adaptive metric. The usage of metric-driven approximation spaces is well-established for the solution of multiscale partial differential equations. Here we demonstrate that it can be also useful for the construction of iterative methods. For the proposed class of conjugate gradient methods we explored the choice of corresponding momentum parameters $\beta^n$ and we numerically investigated five different realizations (including the trivial choice $\beta^n=0$) to compare their performance for the computation of ground states of rotating Bose--Einstein condensates. Among the different options, the Polak--Ribi\'ere and Hestenes--Stiefel realizations performed best and showed the highest efficiency and lowest iteration numbers for the considered application. 



\def\cprime{$'$}


\end{document}